\documentclass[12pt]{article}
\makeatletter\@addtoreset{equation}{section}\makeatother

\newtheorem{Tm}{Theorem}[section]
\newtheorem{Rk}{Remark}[section]
\newtheorem{Lm}{Lemma}[section]
\newtheorem{Co}{Corollary}[section]
\newtheorem{Df}{Definition}[section]
\newtheorem{Pn}{Proposition}[section]

\begin{document}
\title{{{\bf A nonlinear dynamical system on the set of Laguerre entire functions }}}
\author{Yuri Kozitsky$^{a,b,c,1} \ $
 and Lech Wo{\l}owski$^{a,d}$\\ $^{a}$Institute of Mathematics,\\ Marie
Curie-Sk{\l}odowska University, Lublin 20-031 Poland\\ $^{b}$Institute for
Condensed Matter Physics, Lviv 290011 Ukraine \\ $^c$e-mail: {\sf
jkozi@golem.umcs.lublin.pl} \\  $^d$e-mail: {\sf lechw@golem.umcs.lublin.pl}}
\date{}
\maketitle
\par \noindent
{\bf ABSTRACT: } A nonlinear modification of the Cauchy problem $D_t f (t, z)=
\theta D_z f(t,z) + zD^2_z f(t,z)$, $t\in {I\!\! R}_+ = [0, +\infty)$,
 $z\in \hbox{\vrule width 0.6pt height 6pt depth 0pt \hskip -3.5pt}C$,
 $\theta \geq 0 $, $f(0,z) = g(z)\in
{\mathcal L}$ is considered. The set ${\mathcal L}$ consists of Laguerre entire
functions, which one obtains as a closure of the set of polynomials having real
nonpositive zeros only in the topology of uniform convergence on compact subsets
of $ \ \hbox{\vrule width 0.6pt height 6pt depth 0pt \hskip -3.5pt}C$. The
modification means that the time half-line ${I\!\! R}_+$ is divided onto the
intervals ${\mathcal I}_n = [(n-1)\tau , n\tau ]$, $n\in {I\!\! N}$, $\tau>0$,
and on each ${\mathcal I}_n $ the evolution is to be described by the above
equation but at the endpoints the function $f(t, z)$ is changed: $f(n\tau , z)
\rightarrow
\left[ f\left(n\tau , z\delta^{-1-\lambda}\right)\right]^\delta $, with
$\lambda>0$ and an integer $\delta \geq 2$. The resolvent operator of such
problem preserves the set ${\mathcal L}$. It is shown that for $t\rightarrow
+\infty$, the asymptotic properties of $f(t, z)$ change considerably when the
parameter $\tau$ reaches a threshold value $\tau_*$. The limit theorems for
$\tau < \tau_* $ and for $\tau = \tau_* $ are proven. Certain applications,
including limit theorems for weakly and strongly dependent random vectors, are
given.

\begin{large}

\vskip.2cm
\end{large}
%\end{center}
\vskip.1cm
$^{1}${Supported in part under the Grant KBN No 2 P03A 02915}
\section{Setup}
\subsection{Introduction}
The Laguerre entire functions \cite{[Il]} are obtained as uniform limits on
compact subsets of $\ \hbox{\vrule width 0.6pt height 6pt depth 0pt \hskip
-3.5pt}C$
 of the sequences of polynomials
possessing real nonpositive zeros only. These functions are being studied by
many authors during this century in view of their various applications (see also
\cite{[Lev]}). In \cite{[KW]} the set of Laguerre entire functions ${\cal L}$
was described in the framework of locally convex spaces of exponential type
entire functions. In particular, it was shown that the Cauchy problem
\begin{eqnarray*}
\frac{\partial f(t,z)}{\partial t} & = & \theta\frac{\partial f(t,z)}%
{\partial z} +z\frac{\partial^2 f(t, z)}{\partial z^2} , \ \ t\in{I\!\! R}_+
\stackrel{\rm def }{=} [0, +\infty),\ \ z\in \
\hbox{\vrule width 0.6pt height 6pt depth 0pt \hskip -3.5pt}C , \\
f(0,z ) & = & g(z) \in {\cal L}, \ \ \theta \geq 0 .
\end{eqnarray*}
has a unique solution in ${\cal L}$ at least for $t$ small enough. This solution
was obtained in an integral form and its possible asymptotic properties when
$t\rightarrow +\infty$ were considered. In this paper, a nonlinear modification
of this problem is introduced and studied. We divide the time half-line onto the
intervals $[(n-1)\tau, n\tau]$, $n\in {I\!\! N}$ with certain $\tau>0$. On each
such an interval the evolution is to be described by the above equation but at
the endpoints the function $f(t,z)$ is changed
\begin{equation}
\label{i}
f(n\tau, z) \rightarrow \left[f(n\tau, z\delta^{-1-\lambda }) \right]^\delta,
\end{equation}
with a fixed $\lambda >0$ and an integer $\delta \geq 2$. For this dynamical
system, we construct the evolution operator as a holomorphic nonlinear map
between the Fr{\'e}chet spaces of entire functions, which preserves the set of
Laguerre entire functions. Here we use the properties of the operators having
the form $\varphi (\Delta_\theta)$ with $\Delta_\theta = (\theta +z D)D$ and
$\varphi\in {\cal L}$ studied in \cite{[KW]}. For $\lambda <1/2$, we show that,
for sufficiently small values of $\tau$, the asymptotic properties of $f(t,z)$,
$t\rightarrow +\infty$ qualitatively are the same as in the case where the
evolution is described only by the transformation (\ref{i}). At the same time,
it is shown that there exists a threshold value $\tau_* >0$ such that the
asymptotic behaviour of $f(t,z)$ changes drastically when $\tau$ achieves this
value. The description of this phenomenon is based upon the properties of the
evolution operator fixed points. The results obtained are then used to describe
a similar evolution on the sets of isotropic (i.e. $O(N)$--invariant) analytic
functions and measures defined on ${I\!\! R}^N$. In particular, the limit
theorems for strongly and weakly dependent $N$-dimensional random vectors are
proved.

Every statement given below in the form of Proposition either was proved in
\cite{[KW]} or may be proven in an evident way.

\subsection{Definitions and Main Results}

Let ${\cal E}$ be the set of all entire functions $\
\hbox{\vrule width
0.6pt height 6pt depth 0pt \hskip -3.5pt}C \rightarrow
\hbox{\vrule width
0.6pt height 6pt depth 0pt \hskip -3.5pt}C$. For $b>0$, we define
\[
{\cal B}_{b}=\{f\in {{\cal E} }\mid \Vert f\Vert _{b}<\infty \},
\]
where
\begin{equation}  \label{1}
\Vert f\Vert _{b}=\sup_{k\in {I\!\! N}_0 } \{ b^{-k}\mid f^{(k)}(0)\mid \},
\ \ \ f^{(k)}(0)=(D^{k}f)(0) = \frac{d^k f}{dz^k } (0),
\end{equation}
and ${I\!\! N}_0 $ stands for the set of nonnegative integers. For $a\geq 0$,
let
\begin{equation}  \label{5}
{\cal A}_{a}=\bigcap_{b>a}{\cal B}_{b}=\{f\in {\cal E} \mid (\forall b>a ) \
\Vert f\Vert _{b}<\infty \}.
\end{equation}
\begin{Pn}
\label{1pn} $\left( {\cal B}_{b},\Vert \cdot \Vert _{b}\right) $ is a
Banach space, ${\cal A}_a $ equipped with the topology
defined by the family $%
\{\Vert .\Vert _{b}, \ b>a\}$ is a Fr{\'e}chet space.
\end{Pn}
An equivalent topology on ${\cal A}_a$ may be introduced by means of the family
$\{ |. |_b , \ b>a\}$ of the norms
\begin{eqnarray*}
\label{NO2}
|f|_b & \stackrel{\rm def}{=} & \sup_{z\in \ \hbox{\vrule width%
0.6pt height 6pt depth 0pt \hskip -3.5pt}C }%
\{ |f(z)|\exp(-b|z|)\}.
\end{eqnarray*}
\begin{Df}
\label{1df} A family ${\cal L}$  is formed by the entire functions possessing the
representation
\begin{equation}
f(z)=Cz^{m}\exp (\alpha z)\prod_{j=1}^{\infty }(1+\gamma_{j}z),  \label{6}
\end{equation}
\[
C\in \ \hbox{\vrule width 0.6pt height 6pt depth 0pt \hskip -3.5pt}C,\ m\in {
I\!\!N}_{0}, \ \ \alpha \geq 0, \ \  \gamma_{j}\geq \gamma_{j+1}\geq 0 ,\
\sum_{j=1}^{\infty}\gamma_{j}<\infty .
\]
\end{Df}
The elements of ${\cal L}$ are known as the Laguerre entire functions
\cite{[Il]}. Due to Laguerre and P\'olya (see e.g. \cite{[Il]}, \cite{[Lev]}),
we know that ${\cal L}$ consists of the polynomials possessing real nonpositive
zeros only as well as of their uniform limits on compact subsets of $\
\hbox{\vrule width 0.6pt height 6pt depth
0pt \hskip -3.5pt}C$. Let ${\cal P}_{\cal L}$ be the set of polynomials
belonging to ${\cal L}$ and
\begin{equation}
\label{ea1}
{\cal L}^+ \ \stackrel{\rm def }{=} \{f\in {\cal L} \ \mid \ f(0) >0 \}, \ \
{\cal L}^{(1)} \ \stackrel{\rm def }{=} \{f\in {\cal L} \ \mid \ f(0) =1 \},
\end{equation}
\begin{equation}  \label{7}
{\cal L}_{a} \ {\stackrel{{\rm def}}{=}} \ {\cal L}\cap {\cal A}_{a}, \ \ {\cal
L}_{a}^+ \ {\stackrel{{\rm def}}{=}} \ {\cal L}^+ \cap {\cal A}_{a}, \ \ {\cal
L}_{a}^{(1)} \ {\stackrel{{\rm def}}{=}} \ {\cal L}^{(1)}\cap {\cal A}_{a}.
\end{equation}
Given $\theta \geq 0$, a map $\Delta _{\theta }:{\cal E}%
\rightarrow {\cal E}$ is defined to be
\begin{equation}
(\Delta _{\theta }f)(z)=(\theta +zD)Df(z)=\theta \frac{df (z) }{dz}+
z\frac{d^{2}f (z) }{dz^{2}}.  \label{a1}
\end{equation}
For $F(z)=f(z^{2})$, one observes
\begin{equation}
(\Delta _{\theta }f)(z^{2})=\frac{1}{4}\left( \frac{2\theta -1}{z}\frac{dF(z)%
}{dz}+\frac{d^{2}F(z)}{dz^{2}}\right) ,  \label{aa01}
\end{equation}
which means that, for $\theta =N/2$, $N\in {I\!\!N}$, the map (\ref{a1}) is
connected with the radial part of the $N $--dimensional Laplacian
\[
\Delta _{r}=\frac{N-1}{r}\frac{\partial }{\partial r}+\frac{\partial ^{2}}{%
\partial r^{2}}.
\]
Consider now the Cauchy problem:
\begin{eqnarray}  \label{ee1}
\frac{\partial f (t,z)}{\partial t} & = & (\Delta_{\theta}f)(t,z),
 \ \ \ t\in {I\!\! R}_{+} ,
\ z\in \ \hbox{\vrule width 0.6pt height 6pt depth 0pt \hskip -3.5pt}C ,  \\
f(0, z) & = & g(z) ,  \nonumber
\end{eqnarray}
and let the initial condition have the form
\begin{equation}
g(z)=\exp (-\varepsilon z)h(z),\ \ \ h\in {\cal A}_{0},\ \ \varepsilon \geq
0 .  \label{ee3}
\end{equation}
The following statement was proven in \cite{[KW]} as Theorem 1.6.
\begin{Pn}
\label{e1tm}
{\rm (i) \ }  For every $\theta \geq 0$ and $g\in {\cal E}$ having the form
(\ref{ee3}), the problem (\ref{ee1}) has a unique solution in ${\cal
A}_{\varepsilon }$, which possesses the following integral representation
\begin{eqnarray}
\label{ee2}
f(t,z) & = & \exp \left(-\frac{z}{t} \right) \int_{0}^{+\infty }s^{\theta
-1}w_{\theta }\left( \frac{ zs}{t}\right) e^{-s}g(ts)ds,\ \ t>0,
\end{eqnarray}
\begin{equation}
\label{ee4}
w_{\theta} (z) \ \stackrel{\rm def}{=} \sum_{k=0}^{\infty}\frac{z^k }{k!%
 \Gamma (\theta+ k)}.
\end{equation}
\vskip.1cm
\begin{tabular}{ll}
{\rm (ii)} & If in (\ref{ee3}) $\varepsilon >0$,
 the solution (\ref{ee2}) converges in ${\cal A}_\varepsilon $ to zero \\
 & when $t\rightarrow +\infty $.\\[0.1cm]
{\rm (iii)} & If in (\ref{ee3}) $h\in {\cal L}_{0}$ and $\varepsilon =0$,
  the solution (\ref{ee2}) also belongs   \\
& to ${\cal L}_0$. It diverges when $t \rightarrow +\infty$,  which means
$M_f (t, r) \rightarrow +\infty $\\
& for every $r\in {I\!\! R}_+ $. Here
\end{tabular}
\[
M_f (t,r) \ \stackrel{\rm def }{=} \ \sup_{|z|\leq r}|f(t,z)|.
\]
\end{Pn}
By claim (ii), the so called stabilization of solutions holds (see e.g.
\cite{[Kam]} and \cite{[De]}).

We modify the evolution described by the equation (\ref{ee1}) as follows. Let us
divide the time half--line ${I\!\! R}_{+}$ onto the intervals $[(n-1)\tau,
n\tau]$, $n\in {I\!\! N}$ with some $\tau >0$. On each such an interval, the
evolution is to be described by (\ref{ee1}) but at the moments $t=n\tau $,
 $n\in{I\!\! N}_0$ the function is changed as follows
\[
f(n\tau , z) \rightarrow [f(n\tau , z\delta^{-1 -\lambda})]^\delta ,
\]
with a fixed $\lambda > 0$ and an integer $\delta \geq 2$. It is more convenient
to deal with the sequence of functions depending on $t$ from one such interval
instead of considering one function with $t$ varying on the sequence of
intervals. In what follows, we consider the sequence of functions $\{f_n (t,z) ,
n\in {I\!\! N}_0 \}$, each of which is a solution of the following Cauchy
problem
\begin{eqnarray}
\label{ee5}
\frac{\partial f_n (t, z)}{\partial t} & = & \tau (\Delta_\theta f_{n} )(t,z),
\ \ \ \tau \geq 0, \ \ \ t\in [0,1], \ \ z\in \
\hbox{\vrule width 0.6pt height 6pt
depth 0pt \hskip -3.5pt}C  ,\\
f_n (0, z) & = & [f_{n-1} (1, z\delta^{-1-\lambda})]^\delta , \ n\in {I\!\! N}, \nonumber \\
 f_0 (1,z)  & = & g(z) \in {\cal L}^{+} .  \nonumber
\end{eqnarray}
Any $g\in {\cal L}^+ $ is described by the parameters $C$, $\alpha$,
$\{\gamma_j\}$ (see (\ref{6}) and (\ref{ea1})) and one can show that $g\in {\cal
L}_\alpha^+ $. For such functions, we define
\begin{equation}
\label{ee7}
m_k (g) = \sum_{j=1}^{\infty} \gamma_j^k , \ \ \ k\in {I\!\! N},
\end{equation}
and
\begin{equation}
\label{ea7}
I (g) =\left\{
\begin{array}{l}
\lbrack 0,(\delta^\lambda -1) /\alpha],\quad \alpha >0 \\
\\
\lbrack 0,\infty ),\quad \ \alpha =0
\end{array}
\right. .
\end{equation}
Proposition \ref{e1tm} implies the existence of solutions of (\ref{ee5}) at
least for $g\in {\cal L}_0$. The first our theorem establishes the existence of
these solutions for more general situations.
\begin{Tm}
\label{ee0tm}
Let $g\in {\cal L}^+ $ and $\tau \in I(g)$ be chosen. Then for every $n\in
{I\!\! N}$ and $\theta \geq 0$, the problem (\ref{ee5}) has a unique solution
$f_n $, which belongs to ${\cal L}^{+}_\alpha$.
\end{Tm}
 For $\tau =0$, the sequence $\{f_n \}$ can be found explicitly:
\begin{equation}
\label{ea2}
f_n (t,z) = [g(z \delta^{-n(1+\lambda)})]^{\delta^n}.
\end{equation}
If $g\in {\cal L}^{(1)}$, this sequence converges in ${\cal A}_\alpha $ to the
function $f (t, z) \equiv 1$. Thus one may expect that the same or similar
convergence holds also for small positive values of $\tau$. On the other hand,
for large values of $\tau$, claim (iii) of Proposition \ref{e1tm} suggests the
divergence. Our aim in this work is to study the questions: (a) does there exist
the intermediate value of $\tau$, say $\tau_* $, which separates such "small"
and "large" values; (b) what would be the convergence of the sequence $\{f_n \}$
for $\tau = \tau_* $. The answer has been found for the values of $\lambda$
restricted to the interval $\lambda \in (0, 1/2)$ when the initial element $g$
is being chosen in a subset of ${\cal L}^+ $ defined by $\lambda$ as follows.
Let
\begin{equation}
\label{ea6}
\vartheta (\lambda) \ \stackrel{\rm def }{=} \frac{1- \delta^{- \epsilon}}%
{\delta^\lambda -\delta^{-\epsilon}}, \ \  \epsilon = \frac{1-2\lambda}{4} .
\end{equation}
\begin{Df}
\label{2df}
A family ${\cal L}(\lambda)$ consists of the functions $g\in {\cal L}^{(1)} $
which are not constant and are such that
\begin{equation}
\label{ee8}
\frac{m_2 (g)}{[\alpha + m_1 (g) ]^2 } \leq \frac{\delta^{1/2}}{\theta +1}\vartheta%
 (\lambda) , \ \  \ \frac{m_2 (g) }{[m_1 (g) ]^2 } \leq \frac{\delta^{1/2}}{\theta+1} .
\end{equation}
\end{Df}
Thereby, we state our main theorem.
\begin{Tm}
\label{ee1tm}
For every $\theta \geq 0$ and $g\in {\cal L} (\lambda) $, there exist a positive
 $\tau_* \in I(g)$ and a function $C:[0,\tau_* ] \rightarrow {I\!\! R}_+$ such that
\vskip.1cm
\begin{tabular}{ll}
{\rm (i)}  & for $\tau < \tau_*$, the sequence of solutions of (\ref{ee5}) \\
           &  $\{ f_n (t,z) \ \mid \ n\in {I\!\! N}_0 ,  \ f_0 (1, z) = %
              C(\tau)g(z)    \}$ converges \\
           &  in ${\cal A}_{\beta_{*}^{-1} }$,
              $\beta_* \ \stackrel{\rm def}{=} \tau_* / (\delta^\lambda - 1)$
                    to the function $f(t, z) \equiv 1$;\\[.1cm]
          % & \\
{\rm (ii)} & for $\tau = \tau_*$, the sequence
         $\{ f_n (t,z) \ \mid \ n\in {I\!\! N}_0 ,  \ f_0 (1, z) = %
              C(\tau_* )g(z)    \}$  \\
           & converges in ${\cal A}_{\beta_{*}^{-1} }$,
              to
\end{tabular}
\vskip.1cm
\begin{equation}
\label{8}
f_* (t, z) = \delta^{- \delta \theta \lambda / (\delta - 1)}
[1 - t(1-\delta^{-\lambda})]^{-\theta}\exp\left(\frac{1}{\tau_*}
\frac{1-\delta^{-\lambda}}{1-t(1-\delta^{-\lambda})}z \right).
\end{equation}
\end{Tm}
\begin{Rk}
\label{1rk}
The convergence to nontrivial (neither zero nor infinity) limits needs to
control the constant $C$ in the representation (\ref{6}) of the initial element
of $\{f_n
\}$. Otherwise one obtains only such trivial limits for "small" and "large"
values of this constant.
\end{Rk}
\subsection{Some Applications and Further Results}
Let ${\cal E}^{(N)}  $, $N\in {I\!\! N}$ be the set of analytic
functions $F: {I\!\! R}^N \rightarrow \ \hbox{\vrule width 0.6pt%
 height 6pt depth 0pt \hskip -3.5pt}C$. For appropriate
$F\in {\cal E}^{(N)} $ and
some $b>0$, we set
\begin{equation}
\label{9}
\|F\|_{b, N} \ \stackrel{\rm def}{=} \ \sup_{x\in {I\!\! R}^N } \{%
\mid F(x) \mid \exp(-b\mid x\mid^2 ) \},
\end{equation}
where $\mid x \mid $ is the Euclidean norm of $x\in {I\!\! R}^N $.
Let
\begin{equation}
\label{10}
{\cal A}_a^{(N)}  \ \stackrel{\rm def }{=} \ %
\{F\in {\cal E}^{(N)}  \ \mid \ %
\|F\|_{b, N} < \infty, \forall b>a\}, \ \ a\geq 0 .
\end{equation}
This set equipped with the topology generated by the family $\{\| . \|_{b, N},
b>a \}$ becomes a Fr{\'e}chet space. Let $O(N)$ stand for the group of all
orthogonal transformations of ${I\!\! R}^N$. A function $F\in {\cal E}^{(N)} $
is said to be isotropic if for every $U\in O(N)$ and all $x\in {I\!\! R}^N $, $
F(Ux )= F(x)$. The subset of ${\cal E}^{(N )}$ consisting of isotropic functions
is denoted by ${\cal E}_{\rm isot}^{(N )}$. Now let ${\cal P}_{\rm isot}^{ (N
)}\subset%
 {\cal E}_{\rm isot}^{(N )}$ stand for the set of isotropic
polynomials. The classical Study--Weyl theorem \cite{[Weyl]}
(see also \cite{[Lu]}) implies that there exists a bijection between the
set of all polynomials of one complex variable ${\cal P}$ and
${\cal P}_{\rm isot}^{(N )}$ established by
\[
%\label{11}
{\cal P}_{\rm isot}^{(N )}\ni P (x) = p((x,x)) \in {\cal P} ,
\]
where $(. , .)$ is the scalar product in ${I\!\! R}^N$. Obviously each a
function $F$ having the form
\begin{equation}
\label{11x}
F(x) = f((x,x)),
\end{equation}
with certain $f\in {\cal E}$, belongs to ${\cal E}_{\rm isot}^{(N )}$. Given a
subset ${\cal X} \subset {\cal E}$, we write ${\cal X} ({I\!\! R}^N )$ for the
subset of ${\cal E}_{\rm isot}^{(N )}$ consisting of the functions obeying
(\ref{11x}) with $f\in {\cal X}$. In this notation ${\cal P}_{\rm isot}^{(N )} =
{\cal P}({I\!\! R}^N )$. Consider a map
\[
{\cal E}^{(N)}_{\rm isot} \ni F\mapsto \left(\Delta +\frac{d}{(x,x)}(x,\nabla)%
 \right)F\in {\cal E}^{(N)}_{\rm isot},
\]
where $\Delta$ and $\nabla$ stand for the Laplacian and for the gradient in
${I\!\! R}^N $. For a pair of functions $F$ and $f$ satisfying (\ref{11x}), one
has (c.f. (\ref{aa01}))
\begin{equation}
\label{17}
\left(\Delta + \frac{d}{(x,x)}(x, \nabla)\right) F(x)%
 = 4\left(\Delta_{\theta}f\right)((x,x)),
\end{equation}
where $\Delta_\theta $ is defined by (\ref{a1}) with
\begin{equation}
\label{16}
\theta = \frac{ N+ d}{2}.
\end{equation}
Now let us consider the following Cauchy problem -- an analog of (\ref{ee5}):
\begin{eqnarray}
\label{18}
\frac{\partial F_n (t, x)}{\partial t} & = & \tau \left(\Delta + \frac{d}{(x,x)}%
(x, \nabla )\right) F_{n} (t,x),
  \ t\in [0,1],  \ x\in {I\!\! R}^N , \nonumber \\
F_n (0, x) & = & \left[F_{n-1} (1, x\delta^{-(1+\lambda)/2})\right]^\delta %
, \ n\in {I\!\! N},  \\
 F_0 (1,x)  & = & G(x) \in {\cal L}^{+} ({I\!\! R}^N )  .  \nonumber
\end{eqnarray}
For $G\in{\cal L}^{+} ({I\!\! R}^N )$, there exists $g\in {\cal L}^+$ such that
$G$ and $g$ satisfy (\ref{11x}), thus the interval (\ref{ea7}) is defined for
such $G$. The direct corollary of Theorem \ref{ee0tm} reads
\begin{Tm}
\label{2tm}
For every $d\geq -N$, $G\in{\cal L}^{+} ({I\!\! R}^N )$, $\tau \in I(g)$, and
$n\in {I\!\! N}$, the problem (\ref{18}) has a unique solution $F_n $, which
also belongs to ${\cal L}^{+} ({I\!\! R}^N )$.
\end{Tm}
For $\lambda \in (0, 1/2)$, we have an analog of Theorem \ref{ee1tm}.
\begin{Tm}
\label{3tm}
For every $d\geq -N$ and $g\in {\cal L}(\lambda )$, there exist a positive
$\tau_*
\in I(g)$ and $C:[0,\tau_*] \rightarrow {I\!\! R}_+ $, such that
\vskip.1cm
\begin{tabular}{ll}
{\rm (i)}  & for $\tau < \tau_*$, the sequence of solutions of (\ref{18})\\
           &  $\{ F_n (t,x) \ \mid \ n\in {I\!\! N}_0 , \ \ F_0 (1, z) = %
              C(\tau)g((x,x))    \}$ converges  \\
           & in ${\cal A}_{\beta_{*}^{-1} }^{(N)}$,
              $\beta_* \ \stackrel{\rm def}{=} \tau_* / (\delta^\lambda - 1)$
                    to the function $F(t, x) \equiv 1$;\\[.1cm]
          % & \\
{\rm (ii)} & for $\tau = \tau_*$,
         $\{ F_n (t,x) \ \mid \ n\in {I\!\! N}_0 , \ \ F_0 (1, x) = %
              C(\tau_* )g((x,x))    \}$ \\
         &  converges
             in ${\cal A}_{\beta_{*}^{-1} }^{(N)}$
              to
\end{tabular}
\begin{equation}
\label{8x}
F_* (t, x) = \delta^{- \delta \theta \lambda / (\delta - 1)}
[1 - t(1-\delta^{-\lambda})]^{-\theta}\exp\left(\frac{1}{\tau_*}
\frac{1-\delta^{-\lambda}}{1-t(1-\delta^{-\lambda})}(x,x) \right),
\end{equation}
\begin{tabular}{ll}
& where $\theta$ is given by (\ref{16}).
\end{tabular}
\end{Tm}
Let ${\cal M}$ stand for the set of probability measures $\mu$ on ${I\!\! R}^N $
such that
\[
\int_{ {I\!\! R}^N } \exp(\varepsilon (x,x) )\mu(dx ) < \infty ,
\]
with certain $\varepsilon >0$. For each such a measure, the function
\begin{equation}
\label{9y}
F_\mu (x) \ \stackrel{\rm def }{=} \ \int_{{I\!\! R}^N }\exp((x,y))\mu(dy),
\end{equation}
belongs to ${\cal E}^{(N)}$. For a Borel subset $B\subset {I\!\! R}^N$,
 we let
\[
B-x = \{ y\in {I\!\! R}^N \ \mid x+y \in B\}, \ \ UB = \{x\in {I\!\! R}^N%
\ \mid \ U^{-1}x \in B\} , \ U\in O(N) .
\]
A measure $\mu \in {\cal M}$ is said to be isotropic if it is $O(N)$--invariant
(i.e. $\mu(UB) = \mu(B)$), the subset ${\cal M}_{\rm isot}\subset {\cal M}$ is
to consist of such isotropic measures. Obviously,  $F_\mu \in  {\cal
E}^{(N)}_{\rm isot}$ for $\mu
\in {\cal M}_{\rm isot}$. Now let ${\cal M}({I\!\! R}^N  )$ be the subset of
${\cal M}_{\rm isot}$ consisting of the measures for which $F_{\mu} \in  {\cal
L}^{(1)} ({I\!\! R}^N )$. For a pair of measures $\mu, \ \nu \in {\cal M}$,
their convolution is as usual
\begin{equation}
\label{10y}
(\mu \star \nu )(B) = \int_{{I\!\! R}^N }\mu(B-x)\nu(dx).
\end{equation}
Since $F_{\mu \star \nu } = F_\mu F_\nu$, the measure $\mu \star \nu $ belongs
to ${\cal M}({I\!\! R}^N  )$ whenever $\mu$ and $\nu$ possess this property. Now
let $\delta$, $\lambda$, and $\tau$ be as in (\ref{ee5}), (\ref{18}). Consider
the sequence $\{ \mu_n , n\in {I\!\! N}_0 \}$ defined
\begin{equation}
\label{11y}
\mu_n (dy) = \frac{1}{M_n (\tau)}\exp\left(\tau (y,y)\right) \mu_{n-1}^{\star \delta}%
(\delta^{(1+\lambda)/2} dy), \ \ \mu_0 = \nu \in{\cal M}({I\!\! R}^N) ,
\end{equation}
where
\[
M_n (\tau) \ \stackrel{\rm def}{=} \int_{{I\!\! R}^N } %
\exp\left(\tau (y,y)\right) \mu_{n-1}^{\star \delta}%
(\delta^{(1+\lambda)/2} dy),
\]
and $\mu^{\star \delta}$ is the convolution of $\delta$ copies of $\mu$. The
measure $\mu_{n-1}^{\star \delta} (\delta^{(1+\lambda)/2} \cdot)$ describes the
probability distribution of the normalized sum of $\delta$ identically
distributed independent random vectors. By means of the multiplier
$\exp\left(\tau (y,y)\right)$ in (\ref{11y}), we set these vectors being
dependent, thus the measure $\mu_n $ describes the probability distribution of
the following random vector
\begin{equation}
\label{13y}
X^{(n)} = \frac{1}{\sqrt{\delta}} \delta^{-\lambda/2}\left(X^{(n-1)}_1 +\dots +%
X^{(n-1)}_{\delta} \right).
\end{equation}
The normalization of this sum is "abnormal" (more than normal) due to the
additional factor $\delta^{- \lambda/2}$. Every $X^{(m)}$ is the sum of
$\delta^m $ vectors of the zero level. Such random vectors are known to be
hierarchically dependent (see e.g. \cite{[Koz]}). Their dependence is
proportional to the parameter $\tau$ -- it disappears if $\tau
=0$. Therefore, one may expect that, for small positive values of $\tau$, the
dependence remains weak and the vectors obey the classical central limit
theorem. In this case, due to the factors $\delta^{- \lambda/2}$, the sequence
of measures $\{\mu_n \}$ ought to be asymptotically degenerate at zero, which
means that the corresponding by (\ref{9y}) sequence $\{F_{\mu_n }\}$ converges
to the function $F(x) \equiv 1$. But the functions $F_{\mu_n }$ may be obtained
as solutions of the problem (\ref{18}). To use this fact we construct the subset
of ${\cal M}({I\!\! R}^N  )$ corresponding to ${\cal L}(\lambda)$ introduced by
Definition \ref{2df}. Choose $\lambda\in (0,1/2)$. For a measure $\nu \in {\cal
M}({I\!\! R}^N  )$, let $g\in {\cal L}^{(1)}$ be the function such that $F_\nu
(x) = g((x,x))$. Then
\begin{equation}
\label{14y}
{\cal M}_\lambda ({I\!\! R}^N  ) \ \stackrel{\rm def}{ =} \ \{ \nu \in {\cal
M}({I\!\! R}^N )
\ \mid \ g\in {\cal L} (\lambda) \}.
\end{equation}
The following assertion is a corollary of Theorem \ref{3tm} for $d=0$.
\begin{Tm}
\label{4tm}
For every $N\in {I\!\! N}$ and $\nu \in {\cal M}_\lambda ({I\!\! R}^N  )$, there
exists $\tau_* >0$ such that
\vskip.1cm
\begin{tabular}{ll}
{\rm (i)} & for $\tau < \tau_{*}$, the sequence of measures
       defined by (\ref{11y}) \\
& $\{ \mu_n \ \mid \ n\in {I\!\! N}_0 , \ \mu_0 = \nu \}$ converges weakly to the measure \\
& degenerate at zero; \\[.1cm]
{\rm (ii)} & for $\tau = \tau_* $, this sequence converges weakly to the isotropic\\
         &  Gaussian measure with variance $2N(\delta^\lambda -1)
            /\tau_* $.
\end{tabular}
\end{Tm}
It should be pointed out that the convergence to a nondegenerate measure for the
abnormal normalization described by claim (ii) means the appearance of a strong
dependence between the random vectors considered. For $\tau<\tau_*$, the
dependence is weak and the classical central limit theorem ought to hold. To
show this we introduce the classical normalization instead of (\ref{13y}). So we
set along with (\ref{11y}):
\begin{eqnarray}
\label{15y}
\tilde{\mu}_n (dy ) & = &\frac{1}{\tilde{M}_n (\tau )}\exp\left(\delta^{-n\lambda/2}\tau%
(y,y)\right) \tilde{\mu}_{n-1}^{\star \delta}(\sqrt{\delta}dy ), \  \mu_0 = \nu
\in
 {\cal M} ({I\!\! R}^N ), \nonumber \\
\tilde{M}_n (\tau ) & = & \int_{{I\!\! R}^N }\exp\left(\delta^{-n\lambda/2}\tau%
(y,y)\right) \tilde{\mu}_{n-1}^{\star \delta}(\sqrt{\delta}dy ).
\end{eqnarray}
\begin{Tm}
\label{5tm}
Let $N$, $\nu $, and $\tau_* $ be as in Theorem \ref{4tm}. Then for $\tau <
\tau_*$, the sequence of measures $\{ \tilde{\mu}_n \ \mid \ n\in {I\!\! N}_0 ,
\
\tilde{\mu}_0 = \nu\}$
 defined by (\ref{15y}) converges weakly to an isotropic Gaussian measure.
\end{Tm}

\section{Preliminaries}
\subsection{Laguerre Entire Functions and Evolution Operator}

We start with the description of the Fr{\'e}chet spaces ${\cal A}_a$. A subset
$B\subset {\cal A}_{a}$ is said to be bounded in $ {\cal A}_a $ if for every
$b>a$, there exists $K_b >0$ such that, for all $f\in B$, ${\| f \| }_{b}\leq
K_{b}$.
\begin{Pn}
\label{n1tm} For every $a\geq 0$, the space ${\cal A}_{a}$ possesses the
properties: \vskip.1cm
\begin{tabular}{ll}
{\rm (i)} & the relative topology on bounded subsets of $ {\cal A}_{a}$
coincides with \\ & the topology of uniform convergence on compact
 subsets of $\ \hbox{\vrule width 0.6pt height 6pt depth 0pt \hskip -3.5pt}C
$; \\[0.1cm] {\rm (ii)} & multiplication, i.e., $(f,g)\mapsto fg$ is a
continuous map \\ & from ${\cal A}_{a}\times  {\cal A}_{b}$ into ${\cal
A}_{a+b}$; \\[.1cm] {\rm (iii)} & differentiation, i.e. $f\mapsto f' $ is a
continuous self-map of ${\cal A}_a $.
\end{tabular}
\end{Pn}
\begin{Rk}
\label{n1rk}
It can be easily shown that, for positive $a$ and $b$,
\begin{equation}
\label{X1}
\|fg\|_{a+b} \leq \|f\|_a \|g\|_b ,
\end{equation}
thus $(f,g)\mapsto fg$ is a continuous map from
${\cal B}_{a}\times  {\cal B}_{b}$ into ${\cal B}_{a+b}$, which implies
claim {\rm (ii)} of the latter statement.
\end{Rk}
\begin{Pn}
\label{n2tm} Every sequence $\{f_{n},n\in {I\!\!N}\}\subset
{\cal L}_{a}$, $a\geq 0$, that converges in ${\cal E}$ to a function $ f\in
{\cal A}_{a}$, which does not vanish identically, is a bounded subset of ${\cal
A}_{a}$ and hence, by claim (i) of Proposition \ref{n1tm}, it converges in
${\cal A}_{a}$ to $f\in {\cal L}$.
\end{Pn}
For $f\in {{\cal L}}^{+}$, one has $f(0)>0$ (see (\ref{ea1})). Therefore, for
such a function, there exists the neighborhood ${\cal D}$ of the origin in which
$f\neq 0,$ hence the following holomorphic function can be defined
\begin{equation}
\label{19}
\varphi (z)=\log f(z),\qquad z\in {\cal D}.
\end{equation}
In the sequel we use the notation
\begin{equation}
\label{deriv}
\varphi^{(k)}=(D^{k}\log f)(0),\qquad k\in {I\!\! N}_0 .
\end{equation}
\begin{Pn}
{\rm [The sign rule]} \label{SRule} Let $f\in {\cal L}^{+},$ then
\begin{equation}
\label{fisr}
(-1)^{k-1}\varphi ^{(k)}\geq 0,\quad k\in {I\!\!N} .
\end{equation}
Equalities hold simultaneously for all $k\geq 2$%
\ and only for $f(z)= C \exp (az).$
\end{Pn}
\begin{Lm}
\label{Conver}
For a sequence $\{f_{n}(z)\ | \ n\in%
{I\!\! N} _{0}\ \ f_{n}\in {\cal L}^{+}\},$\ let the derivatives (\ref{deriv})
 satisfy:{\rm (i)  }  $\{ \varphi_n^{(k)} \}$ converges to $\varphi^{(k)}$, $k=0,1$; {\rm
(ii) } $\{\varphi _{n}^{(2)}\}$ converges to zero. Then $\{f_{n}\}$\ converges
to $\exp (\varphi^{(0)} + \varphi^{(1)} z)$\ in ${\cal A}_{\psi}$,  $\psi =
\sup
\varphi_n^{(1)}$.
\end{Lm}

{\bf Proof.} By claim (ii) of Proposition \ref{n1tm} and Proposition \ref{n2tm},
to prove this statement we only need to show that the sequence $\{ f_n (z) / f_n
(0) \}$ converges to $\exp(\varphi^{(1)} z)$ uniformly on compact subsets of $ \
\hbox{\vrule width 0.6pt height 6pt depth 0pt \hskip -3.5pt}C$. Due to known
Vitali's theorem and to the fact that, for the functions considered, $M_f (r) =
f(r)$, we may do this by proving the pointwise convergence of $%
\{f_{n} (z) / f_n (0) \}$ on ${I\!\!R}_{+}$.
To this end we use the specific form of $f\in {\cal L}^{(1)}$ given by
(\ref{6}). For each $\gamma \geq 0 $, one has $\exp (\gamma -\frac{1}{2}\gamma
^{2})\leq 1+\gamma \leq \exp (\gamma )$. Hence for $z\in {I\!\!R}_{+},$
\begin{equation}
\label{20}
\exp (z{\varphi }_{n}^{(1)}+\frac{1}{2}z^{2}{\varphi }_{n}^{(2)})\leq
\frac{f_{n}(z)}{f_n (0)}\leq \exp (z{\varphi }_{n}^{(1)}),
\end{equation}
 which yields the stated convergence.\qquad \hfil\vrule height6pt
width5pt depth-1pt

For an entire function $f\in {\cal E}$ and $t \geq 0$, we define
\begin{equation}
\label{21}
(\exp(t \Delta_{\theta} )f) (z) = \sum_{k=0}^{\infty}\frac{t^k }{k!}%
(\Delta_{\theta}^k f)(z).
\end{equation}
\begin{Pn}
\label{X1pn}
For every positive $a$ and $t$ obeying $a t<1 $, and $\theta \geq 0$,
\[
\|\exp(t \Delta_{\theta} )f \|_b \leq (1- a t)^{-\theta} \|f\|_a, \ \ %
b=a/(1-a t) ,
\]
which means that (\ref{21}) defines a continuous linear map
\begin{equation}
\label{XX1}
{\cal A}_a \ni f \mapsto f_t \ \stackrel{\rm def }{=} \ (\exp(t
\Delta_{\theta} )f ) \in {\cal A}_b , \ \ b = a/(1-at).
\end{equation}
\end{Pn}
\begin{Co}
\label{XX1co}
For every positive $a$ and $t_0 $, a map $ (0,t_0 )\ni t \mapsto f_t \in  {\cal
A}_{b_0 }$, where $b_0  \stackrel{\rm def}{=} a/(1-at_0 )$ and $ f_t $ is
defined by (\ref{XX1}), is differentiable on $(0,t_0 )$ and
\begin{equation}
\label{XX2}
\frac{\partial f_t}{\partial t} = \Delta_\theta f_t , \ \ \ t\in (0,t_0 ).
\end{equation}
\end{Co}
One of the main results of \cite{[KW]} is Theorem 1.3  which asserts that the
operators of the type of (\ref{21}) preserves the class ${\cal L}$. In our case
it reads as follows
\begin{Pn}
\label{X2pn}
Let $a$, $b$, $t$, and $\theta$ be as in Proposition \ref{X1pn}. Then the
operator (\ref{21}), (\ref{XX1}) maps ${\cal L}_a $ into ${\cal L}_b $.
\end{Pn}
The following statements have also been proven in \cite{[KW]}.
\begin{Pn}
\label{X3pn}
For $t>0$, the above operator has the integral form:
\begin{eqnarray}
\label{22}
(\exp(t \Delta_{\theta} )f) (z) =  \exp \left(-\frac{z}{t} \right)
 \int_{0}^{+\infty}s^{\theta -1}w_\theta \left(\frac{zs}{t}\right)e^{-s}%
f(t s)ds ,
\end{eqnarray}
where $w_\theta$ is defined by (\ref{ee4}).
\end{Pn}
\begin{Rk}
\label{n2rk}
Setting in (\ref{22}) $z=0$, one obtains for $f\in {\cal L}^{+}$ and $\theta>0
$: $ \ (\exp(t \Delta_{\theta})f)(0) >0$. On the other hand, one has from
(\ref{21})
\begin{eqnarray}
\label{aa}
 (\exp(t \Delta_{\theta} )f) (0)=\sum_{k=0}^{\infty}\frac{t^k }{k!}f^{(k)}(0)
\frac{\Gamma (\theta + k)}{\Gamma (\theta)}.
\end{eqnarray}
Passing here to the limit $\theta \rightarrow 0$ one gets
\begin{eqnarray}
\label{ab}
 (\exp(t \Delta_{0} )f) (0)= f(0)>0.
\end{eqnarray}
Below the case $\theta =0$ is always understood as the above limit.
\end{Rk}
\begin{Pn}
\label{Shift}
Let$\ v\in {I\!\! R}$ and $\exp (vz)h(z)\in {\cal A}_{b}$, ( $b\geq 0$). For any
$u\geq 0$ obeying the condition $u b<1$,
\begin{equation}
\label{ident}
\exp (u{\Delta }_{\theta })\exp (vz)h(z)=\exp \left(\frac{vz}{1-uv }\right)h_u (z),
\end{equation}
where
\begin{equation}
\label{id2}
h_u ( z)=(1-uv)^{-\theta }\exp (u(1-uv){\Delta }_{\theta })h\left(\frac{z}{%
(1-uv)^{2}}\right).
\end{equation}
Moreover, if $h\in { \cal A}_{a}$, then $h_u \in {\cal A}%
_{c},$ where
\begin{equation}
\label{cident}
c=a(1-uv)^{-1}(1-(v+a)u)^{-1}.
\end{equation}
\end{Pn}
By means of (\ref{21}), we construct the evolution operator which solves
(\ref{ee5}):
\begin{equation}
\label{23}
f_n (t,z) = \exp\left(t\tau \Delta_\theta \right)\left[f_{n-1}(1, z\delta^{-1-\lambda})%
\right]^\delta \ \stackrel{\rm def}{=} T_{t} (f_{n-1} (t, \cdot )) (z) ,
\end{equation}
provided all $f_n (t, z) $ belong to the domain of $T_t $, $t\in [0,1]$. For
short we write
\begin{equation}
\label{24}
f_n (1,z) \ \stackrel{\rm def}{=} \ f_n (z) , \ \  T_1 \ \stackrel{\rm def}{=} \ T.
\end{equation}
Thus one has
\begin{equation}
\label{X5}
f_n = T(f_{n-1 }).
\end{equation}
We use such $T_t $ to define the operators between the Fr{\'e}chet spaces ${\cal
A}_a $ and the Banach spaces ${\cal B}_b$. In all such cases we denote them by
$T_t $ pointing out if necessary between which spaces acts given $T_t $.
Combining claim (ii) of Proposition \ref{n1tm}  with Propositions \ref{X1pn} and
\ref{X2pn} one has
\begin{Pn}
\label{5Xpn}
For every $a<\delta^\lambda /t\tau $, the operator $T_t $ continuously maps:
${\cal B}_a \rightarrow {\cal B}_b \ $,${\cal A}_a \rightarrow {\cal A}_b $, and
${\cal L}_a \rightarrow {\cal L}_b $, where
$b = a\delta^{-\lambda}/ (1-a t\tau\delta^{-\lambda})$.
\end{Pn}
\begin{Pn}
\label{X4pn}
Let the sequence $\{f_n (z) \ \mid \ n\in {I\!\! N}_0 , \  f_0 (z) = g(z) \in
{\cal L}^+ \}$ defined by (\ref{X5}) converge in ${\cal A}_a$, $a\geq 0$ to a
function $f$. Then the sequence of solutions of (\ref{ee5}) $\{ f_n (t, z) \ | \
n\in {I\!\! N}_0 , \ f_0 (z) = g(z)\}$, defined by (\ref{23}), converges in
${\cal A}_a $ to the function
\begin{equation}
\label{X4}
f(t, z) = (T_t f)(z).
\end{equation}
\end{Pn}
To establish the existence and convergence of $\{ f_n \}$ we use an analog of
the Fr{\'e}chet derivative of $T$ and then study the fixed points of $T$ and
their stability. The following corollary of Proposition \ref{5Xpn} allows to
define the differentiability of $T$ acting between the Fr{\'e}chet spaces. For
$a\in [0, \tau
\delta^{-\lambda})$, we set
\begin{equation}
\label{bofa}
b(a) \ \stackrel{\rm def }{=} \frac{a\delta^{-\lambda}}{1- a\tau \delta^{-\lambda}}.
\end{equation}
\begin{Co}
\label{5Xco}
Let $a<\delta^\lambda /\tau $, then there exists $\varepsilon >0$ such that, for
all $a' \in (a, a+\varepsilon)$, the operator $T$ continuously maps ${\cal
B}_{a'}$ into ${\cal B}_{b(a')}$.
\end{Co}
\begin{Df}
\label{5Xdf}
The operator $T : {\cal A}_a \rightarrow {\cal A}_{b(a)} $
is said to be differentiable on ${\cal A}_a $ if for every $f\in {\cal A}_a$,
there exist $\varepsilon >0$ and
a continuous linear operator $T' [f]:{\cal A}_a \rightarrow {\cal A}_{b(a)} $
such that, for every $a' \in(a, a+\varepsilon)$, $T' [f]$ is the Fr{\'e}chet derivative
of $T $ considered as an operator between the Banach spaces
${\cal B}_{a'}$ and ${\cal B}_{b(a')}$. The operator $T' [f]$ is said to be a derivative
of $T$ at $f$.
\end{Df}
\begin{Lm}
\label{DerT}
For $a<\delta^\lambda /\tau $,  the operator $T : {\cal A}_a \rightarrow {\cal
A}_{b(a)} $ is differentiable on ${\cal A}_{a}$\ and its derivative $T^{\prime
}[f]$ is the following operator
\begin{equation}
\label{F'}
T^{\prime }[f](h)(z) =\delta  \exp%
(\tau\Delta _{\theta })\left( (f^{\delta -1}h)(\delta^{-1-\lambda} z )\right).
\end{equation}
\end{Lm}

{\bf Proof.} For $a^{\prime }\in (a,  \delta^\lambda /\tau )$ and $f,$ $h\in
{\cal B}_{a^{\prime }},$ one has
\[
T(f+h)= T(f) + \delta \exp (\tau\Delta%
_{\theta })\left( (f^{\delta -1}h)(\delta^{-1-\lambda} \cdot )\right) + R(f, h),
\]
\begin{eqnarray*}
R(f,h) & = &\exp (\tau\Delta _{\theta }) \left(%
\sum_{k=2}^{\delta }{{{\delta }}\choose{{k}}}f^{\delta -k}h^{k}\right)
(\delta^{-1-\lambda} \cdot ).
\end{eqnarray*}
By means of Remark \ref{n1rk}, (\ref{X1}), and Proposition \ref{X1pn}, one
obtains
\begin{eqnarray*}
\left\|\exp(\tau\Delta_{\theta})\left(f^{\delta- k}h^k \right)(\delta^{-1-\lambda}\cdot )%
\right\|_{b(a' )} & \leq & (1- a' \tau\delta^{-\lambda})^{-\theta}%
\|f\|_{a' }^{\delta -k}\|h\|_{a' }^k , \\
 & &  k= 1, 2, \dots , \delta .
\end{eqnarray*}
This gives for all $a' \in (a, \delta^\lambda / \tau )$,
\[
\|R(f,g )\|_{b(a' )} = o (\|h\|_{a' }),
\]
and also for $T'$ defined by (\ref{F'}),
\[
\|T' [f] (h) \|_{b(a' )} \leq \delta (1- a' \tau\delta^{-\lambda})^{-\theta}%
 \|f\|_{a' }^{\delta - 1} \|h\|_{a' }.
\]
By the latter estimate, $T' [f]$ continuously maps ${\cal B}_{a' }$ into ${\cal
B}_{b(a' )}$ whereas by the former one, this operator is the Fr{\'e}chet
derivative of $T : {\cal B}_{a' }\rightarrow {\cal B}_{b(a' ) }$.
\qquad \hfil\vrule height6pt
width5pt depth-1pt
\newline
The case of $\tau = 0$ was considered in (\ref{ea2}), thus from now on we assume
$\tau >0$. It turns out that it is convenient to remove the explicit dependence
on $\tau$ from the operator $T$. To this end we set
\begin{equation}
\label{Z1}
\tau \ \stackrel{\rm def }{=} \beta (\delta^\lambda -1),
\end{equation}
and include $\beta$ into $z$. Then we consider the sequence $\{g_n (z) \}$
\begin{eqnarray}
\label{Z2}
g_n (z) & =  & Q (g_{n-1})(z) , \ \ \ \ n\in {I\!\! N}, \\ &\stackrel{\rm
def}{=} &
\exp\left((\delta^\lambda -1)\Delta_\theta \right)\left[g_{n-1} (\delta^{-1
-\lambda} z)\right]^\delta , \ g_0
(z) = g(\beta z) , \nonumber
\end{eqnarray}
where $g$ is the function which defines the starting element of $\{f_n \}$. To
prove Theorem \ref{5tm} we shall also use the sequence of functions from ${\cal
L}^{(1)}$,
 $\{ \tilde{g}_n (z) \ | \ n\in {I\!\! N}_0 , %
\ \tilde{g}_0 (z) = g(\beta z)\}$, where $g$ is as above, and
\begin{eqnarray}
\label{Z7a}
\tilde{g}_n (z) & = & \tilde{Q}_n (\tilde{g}_n )(z) \\
& \stackrel{\rm def}{=} & \frac{1}{\tilde{Y}_n }\left\{
\exp\left(\delta^{-n\lambda } (\delta^\lambda -1)\Delta_\theta
\right)\left[\tilde{g}_{n-1}(\delta^{-1}\cdot )\right]^\delta \right\}(z),
\nonumber \\
\tilde{Y}_n &\stackrel{\rm def}{=} & \left\{ \exp\left(\delta^{-n\lambda }%
(\delta^\lambda - 1)\Delta_\theta \right)\left[\tilde{g}_{n-1}(\delta^{-1}\cdot
)\right]^\delta
\right\}(0).
\nonumber
\end{eqnarray}
Comparing (\ref {23}), (\ref{24}) with (\ref{Z2}) one obtains from Proposition
\ref{5Xpn} and Lemma \ref{DerT}.
\begin{Pn}
\label{newpn}
For every $a< \delta^\lambda /(\delta^\lambda -1)$, $Q$ is a differentiable (and hence
continuous) operator, which maps: ${\cal A}_a \rightarrow {\cal A}_{b'}$,
${\cal L}_a^{+} \rightarrow {\cal L}_{b'}^{+}$, where $b' = a [\delta^\lambda - a(%
\delta^\lambda -1)]^{-1}$. Its derivative is
\begin{eqnarray}
\label{Z3}
Q'[g] (h)(z) & =  & \exp\left((\delta^\lambda -1)\Delta_\theta \right)%
\left(\left[g^{\delta - 1}h \right](\delta^{-1-\lambda}z)\right) .
\end{eqnarray}
\end{Pn}
For $\tau \in I(g)$, $\beta$ varies in $J(g) \stackrel{\rm def }{=}%
 (0, 1/\alpha ]$ (see (\ref{Z1}) and (\ref{ea7})).
Let $g\in {\cal L}^+ $ be chosen. Then it possesses the representation (\ref{6})
with $\alpha \geq 0 $, which determines the intervals $I(g)$ (\ref{ea7}) and
$J(g)$, and $g\in {\cal L}_\alpha^+%
\subset{\cal A}_\alpha $.
\begin{Lm}
\label{Y1lm}
For $\tau \in I(g)$, all the elements of $\{f_n \ | \ n\in {I\!\! N}_0 ,
\ f_0 = g \}$ belong to ${\cal L}_\alpha^+\subset{\cal A}_\alpha $,
all the elements of $\{g_n \ | \ n\in {I\!\! N}_0 , \ g_0 (z) = g(\beta z) \}$ belong to ${\cal
L}^{+}_{\beta \alpha}$.
\end{Lm}

{\bf Proof.} For $\tau \in I(g)$, $\alpha \leq (\delta^\lambda
- 1)/\tau <%
\delta^\lambda /\tau$, thus by Corollary \ref{5Xco}, $T$ maps ${\cal A}_\alpha$ into
${\cal A}_{b(\alpha )}$ with
\[
b(\alpha ) = \frac{\alpha \delta^{-\lambda} }{1 - \alpha \tau \delta^{-\lambda }}%
\leq \frac{\alpha \delta^{-\lambda} }{1-1 + \delta^{-\lambda}} = \alpha ,
\]
which means $T:{\cal A}_\alpha \rightarrow {\cal A}_\alpha $. By Proposition
\ref{5Xpn}, $T$ maps ${\cal L}$ into itself; by Remark \ref{n2rk}, $(T f)(0)
>0 $ whenever $f(0) >0$. The second part of Lemma concerning $\{ g_n \}$
directly follows from the first one.
\qquad \hfil\vrule height6pt
width5pt depth-1pt
\newline Since the starting element of $\{g_n \}$ is of
the form $g_0 (z) = g(\beta z)$, all its elements depend on $\beta$. Therefore,
one may consider $g_n $ as a map from $J(g)$ into ${\cal A}_1$. To emphasize
this fact we write sometimes $g_n (\cdot , \beta)$ instead of $g_n $.
\begin{Lm}
\label{DergB} For every $n\in {I\!\! N}_0 $, the map
\begin{equation}
g_{n}:J(g)\rightarrow {\cal A}_{1}
\label{gIwA}
\end{equation}
is differentiable on ${\rm Int}J(g)$, its derivative at $\beta $ is an entire
function $\dot{g}_{n}\in{\cal A}_1 $.
\end{Lm}

{\bf Proof.} Let us show that, for $\beta \in {\rm Int}J(g)$, there exists an
entire function $\dot{g}_{n}\in {\cal A}_1 $ such that, for $\tilde{\beta}
\in {\rm Int}J(g)$,%
\begin{equation}
\label{defdiff}
g_{n}(\cdot , \tilde{\beta})-g_{n}(\cdot ,\beta )=\Delta \beta \dot{g}_n %
+r_{n}(\cdot ,\Delta \beta ), \qquad  \Delta \beta =\tilde{\beta}-\beta ,
\end{equation}
where $r_{n}(\cdot ,\Delta \beta)/\Delta \beta \rightarrow 0$ in ${%
\cal A}_{1} $ when $\Delta \beta \rightarrow 0$.
By claim (iii) of Proposition \ref{n1tm}, differentiation is a continuous
self-map of ${\cal A}_a$. Since $g_0 (z,\beta ) = g(\beta z)$, the stated
property obviously holds for $n=0$. For some $n\geq 1$, let $\dot{ g}_{n-1}$
obey (\ref{defdiff}) and belong to ${\cal A}_{1}$. Then
\begin{eqnarray}
\label{Tg1}
& & g_{n}(\cdot ,\tilde{\beta})-g_{n}(\cdot ,\beta )  =  Q(g_{n-1} (\cdot ,\tilde{\beta} ))
-Q(g_{n-1}(\cdot ,\beta ) )\\
& = & Q\left[(g_{n-1})(\cdot ,\beta) +\Delta \beta \dot{g}_{n-1} +r_{n-1}(\cdot
,\Delta \beta )\right]
-Q(g_{n-1}(\cdot ,\beta)). \nonumber
\end{eqnarray}
By means of the derivative (\ref{Z3}), it can be written as
\begin{eqnarray*}
g_{n}(\cdot ,\tilde{\beta})-g_{n}(\cdot ,\beta ) & = & \Delta \beta Q' [g_{n-1}]\left(%
\dot{g}_{n-1}\right) + Q' [g_{n-1}]\left(r_{n-1}(\cdot ,\Delta \beta )\right) +
R_{n-1}, \nonumber
\end{eqnarray*}
where for all $a>1$,
\begin{eqnarray*}
\|R_{n-1}\|_a & = & o\left(\Delta\beta \|\dot{g}_{n-1}\|_{c(a)} +
\|r_{n-1}(\cdot, \Delta\beta )\|_{c(a)} \right) \\
& = & o(\Delta\beta ), \ \ \ \ \ c(a) \ \stackrel{\rm def}{=}
\frac{a\delta^\lambda} {1+ a(\delta^\lambda -1)}.
\end{eqnarray*}
Since the operator $Q'[g_{n-1}]$ is linear and continuous, the function
\[
 Q' [g_{n-1}]\left(r_{n-1}(\cdot ,\Delta \beta)\right) +
R_{n-1}
\]
obeys the conditions imposed on $r_n$, thus $\dot{g}_n $ exists and
\begin{equation}
\label{Z5}
\dot{g}_n = Q' [g_{n-1}]\left(%
\dot{g}_{n-1}\right).
\end{equation}
\qquad \hfil\vrule height6pt
width5pt depth-1pt
\newline Let $g_n^{(k)} \stackrel{\rm def }{=} D^k_z g_n $, $k\in{I\!\! N}$,
then claim (iii) of Proposition \ref{n1tm} implies
\begin{Co}
\label{DiffB}
For every $n\in {I\!\! N}_0 $ and $k\in {I\!\! N}$, the map
$g_{n}^{(k)}:J(g)\rightarrow {\cal A}_{1}$ is differentiable on ${\rm Int}J(g)$,
its derivative at $\beta $ is an entire function $\dot{g}_{n}^{(k)}$ from ${\cal
A}_1 $. For every $z_0 \in
\hbox{\vrule width 0.6pt
height 6pt depth 0pt \hskip -3.5pt}C$,
$g_n^{(k)}(z_0 , \beta)$ is $ \beta$--differentiable on ${\rm Int}J(g)$ and
\begin{equation}
\label{Z6}
\frac{\partial g_n ^{(k)} (z_0 ,\beta)}{\partial \beta} = \dot{g}_n^{(k)} (z_0 , \beta).
\end{equation}
\end{Co}

\subsection{Invariant Sets and Fixed Points}
By Lemma \ref{Y1lm}, for chosen $g\in {\cal L}^+$ and $\tau
\in I(g)$,  ${\cal L}_\alpha^+$ is an invariant set of $T$. It turns out that
this set contains a subset which $T$ maps into itself as well. Proposition
\ref{Shift} implies that such one is
\begin{equation}
\label{G1}
{\cal G} \ \stackrel{\rm def}{=} \ \{ f(z) = C\exp (uz) \ | \ C>0 ,
 \ u\geq 0 \}\subset {\cal L}^+ .
\end{equation}
In fact
\begin{equation}
\label{G2}
T\left(C\exp (uz) \right) = C^\delta (1- u\tau \delta^{-\lambda} )^{-\theta}%
\exp\left( \frac{u\delta^{-\lambda}z}{1- u\tau \delta^{-\lambda} } \right),
\end{equation}
which also yields that ${\cal G}$ contains the following fixed points of $T$:
\begin{eqnarray}
\label{Y2}
& & f_{i,*} (z)  =  C_{i,*} \exp (u_{i,*} z), \ \ i=1,2, \\
\label{Y2a}
& &  C_{1,*}  =  1,  \ u_{1,*} = 0; \ \ C_{2,*} =%
\delta^{-\lambda\theta /(\delta -1)}, \ u_{2,*} = \frac{1}{\tau}(\delta^\lambda -1).
\end{eqnarray}
Consider the sequence $\{ f_n \ | \ n\in {I\!\! N}_0 , \ f_0 = C_0 g = C_0
\exp(\alpha z) \in {\cal G} \}$. By means of (\ref{G2}), one can calculate $f_n$
explicitly
\begin{eqnarray}
\label{G3}
&  & f_n (z)  =  C_n \exp(u_n z) , \\
&  & C_n  =  C_0^{\delta^n}\Xi_n , \ \ \ \ \ u_n = \frac{\alpha \delta^{-n\lambda}}%
{1-\frac{\alpha \tau }{\delta^\lambda - 1}(1-\delta^{-n\lambda})}, \nonumber \\
&  & \Xi_n  =  \xi_n \prod_{l=1}^{n-1} \xi_l^{(\delta - 1)\delta^{n-1-l}}, \
\xi_l = \left[1-\frac{\alpha \tau }{\delta^\lambda - 1}(1-\delta^{-l\lambda})%
\right]^{- \theta}. \nonumber
\end{eqnarray}
In this case we may check the validity of Theorem \ref{ee1tm} directly. Set
\begin{equation}
\label{GG1}
\tau_* \ \stackrel{\rm def}{=} \frac{1}{\alpha}(\delta^\lambda -1),
\end{equation}
\begin{equation}
\label{GG2}
C(\tau) \ \stackrel{\rm def}{=} \ \prod_{k=0}^{\infty}\left(\frac{\delta^\lambda
-1 -%
\alpha \tau + \alpha \tau \delta^{- (k-1)\lambda}}{\delta^\lambda -1 -%
\alpha \tau + \alpha \tau \delta^{- k\lambda}}\right)^{\theta\delta^{-k-1}} .
\end{equation}
Then for $\tau <\tau_* $, the sequence $\{f_n \ | \ n\in {I\!\! N}_0 , f_0 (z)%
= C(\tau)\exp(\alpha z) \}$ converges in ${\cal A}_\alpha$ to $f_{1,*} \equiv 1$.
If for such $\tau$, one chooses $f_0 (z) = C_0 \exp(\alpha z)$ with $C_0 <
C(\tau)$ (resp. $C_0 > C(\tau $)), then $C_n $ in (\ref{G3}) tends to zero
(resp. to infinity). For $\tau = \tau_* $, one has in (\ref{GG2}) and (\ref{G3})
respectively
\begin{eqnarray*}
&  & C(\tau_* )  =  C_{2,*}, \\
&  & C_n   =  C_0^{\delta^n }\exp\left(\lambda\theta%
 \frac{\delta^n - 1}{\delta -1}\log \delta\right), \ \ u_n = \alpha .
\end{eqnarray*}
Thus for all $n\in {I\!\! N}_0 $, $C_n = C_{2,*}$ if $C_0 = C(\tau_* ) =
C_{2,*}$. For $C_0 < C_{2,*}$ (resp. $C_0 > C_{2,*}$), $C_n $ tends to zero
(resp. to infinity). The fixed points of $Q$ in ${\cal G}$ are
\begin{equation}
\label{Y3a}
g_{i,*} (z) = C_{i,*} \exp(v_{i,*} z), \ \ v_{1,*} = 0 , \ \ v_{2,*} = 1 .
\end{equation}
To describe the stability of the fixed points (\ref{Y2}) we solve the eigenvalue
problem
\begin{equation}
\label{Y4}
T' [f_{i,*}](h) = \Lambda h.
\end{equation}
To this end we set
\[
h(z) = f_{i,*}(z) p(z) = C_{i,*} \exp(u_{i,*} z) p(z),
\]
with $p$ being a polynomial, and obtain from (\ref{F'}) and Proposition \ref{Shift}
\begin{eqnarray*}
T' [f_{i,*}](h)  & =  & \delta C_{i,*}^\delta (1- u_{i,*} \delta^{-\lambda}\tau )^{-\theta}%
\exp\left(\frac{u_{i,*} \delta^{-\lambda}z}{1- u_{i,*} \delta^{-\lambda} \tau } \right) \\
& & \exp\left( \tau(1-u_{i,*} \delta^{-\lambda}\tau)\Delta_\theta \right)p\left(\frac{z\delta^{-1-\lambda}}%
{(1-u_{i,* }\delta^{-\lambda}\tau)^2}\right).
\end{eqnarray*}
Suppose that $\deg p = k$,
$k\in {I\!\! N}_0$ and apply the latter formula in (\ref{Y4}).
Since $\exp(\dots \Delta_\theta )$ maps such $p$ into a polynomial $q$, $\deg q = k$, we may find
$\Lambda^{(i)}_k $ setting the coefficients in front of the $k$-th powers of $z$ to be equal.
This yields
\begin{equation}
\label{Y5}
\Lambda^{(i)}_k = \frac{\delta^{-k\lambda -k +1}}{(1-u_{i,*} \delta^{-\lambda}\tau)^{2k}} ,\ \
 k\in {I\!\! N}_0 .
\end{equation}
For both $ f_{i,*}$, $\Lambda_0 = \delta >1$, which corresponds to their
instability with respect to the variation of the constant multiplier $C$. The
rest of the eigenvalues which describe $f_{1*}$ are $\Lambda^{(1)}_k =
\delta^{-k\lambda -k+1} <1$. But for $f_{2,*}$, one has
\begin{equation}
\label{Y6}
\Lambda^{(2)}_k = \delta^{k\lambda -k+1}, \ \ k\in {I\!\! N}_0 .
\end{equation}
The eigenvalues of $ Q' [g_{i,*}]$ are exactly the same as given by (\ref{Y5}).
For $\lambda \in (0,1/2)$, $\Lambda^{(2)}_2 <1$. This means that, in the
corresponding spaces ${\cal A}_{a}$, $f_{2,*}$, $g_{2,*}$ have the stable
manifolds of ${\rm codim} =2$ and $f_{1,*}$, $g_{1,*}$ have those of ${\rm
codim}
=1 $. This fact plays an important role in proving the convergence to these
fixed points. The proof will be done by showing that there exist $\beta_* >0$
and a function $C: (0,\beta_* ]
\rightarrow {I\!\! R}_+  $ such that all elements of the sequence $\{ g_n \ | \
g_0 (z) = C (\beta_* ) g(\beta_* z)\}$ remain in the stable manifold of
$g_{2,*}$ and the elements of $\{ g_n \ | \ g_0 (z) = C(\beta) g(\beta z), \
\beta <\beta_* \}$ remain in the stable manifold of $g_{1,*}$. The convergence
of the corresponding sequences $\{f_n \}$ will be obtained as a direct
corollary.
\section{Proofs}

\subsection{Main Lemmas}
The case where the initials elements of the sequences considered are chosen in
${\cal G}$ has already been described, thus from now on we suppose that these
functions are chosen outside of ${\cal G}$. We set (see (\ref{19}),
(\ref{deriv}))
\begin{equation}
\label{Z7b}
g_n (z) = C_n \exp(\varphi_n (z) ),  \ \ \tilde{g}_n (z) =
\exp(\tilde{\varphi}_n (z) ),
\ \ \varphi_n (0) = \tilde{\varphi}_n (0) = 0,
\end{equation}
and for $ k\in {I\!\! N}$,
\begin{eqnarray}
\label{Z8}
\varphi_{n}^{(k)} \ \stackrel{\rm def}{=} \  (D^{k}{\varphi}_n )(0), \ \
\tilde{\varphi}_{n}^{(k)} \ \stackrel{\rm def}{=} \  (D^{k} %
\tilde{\varphi}_{n})(0) = \delta^{n\lambda k} \varphi_{n}^{(k)}.
\end{eqnarray}
As it has been shown above (Lemma \ref{DergB} and Corollary \ref{DiffB}), all
$\varphi_n^{(1)}$ are differentiable, and hence continuous, functions of $\beta
\in J(g)$. For $\beta=0$, all $\varphi_n^{(1)}$ are equal to zero, thus one may
assume that, for every $n\in {I\!\! N}_0$, the following inequality
\begin{equation}
\label{Ineq}
\varphi_n^{(1)} <(1-\delta^{-\lambda})^{-1},
\end{equation}
holds for $\beta$ small enough, say, for $\beta\in J_n = (0, \hat{b}_n) $. Below
we obtain the estimates which allow to evaluate the intervals $J_n $. Thus we
set
\begin{equation}
\label{Z9}
{\nu }_{n}=\frac{1}{1-(1-\delta ^{-\lambda }){\varphi}_{n-1}^{(1)}}%
;\quad {\kappa }_{n}=\delta ^{-\lambda }{\nu }_{n}.  \label{kap}
\end{equation}
\begin{Lm}{\rm [Main estimates]}
\label{Z2lm}
For  $\beta \in J_n $, the following estimates hold
\begin{eqnarray}
\label{ce1}
C_n & \geq & C_{n-1}^\delta; \\
\label{ce2}
C_n & \leq &  C_{n-1}^\delta\nu_n^\theta ;
\end{eqnarray}
\begin{equation}
\label{ce3}
C_n \geq C_{n-1}^\delta\nu_n^\theta \exp\left\{ \frac{1}{2} \theta (\theta +1)%
\kappa_n^2 (1- \delta^{- \lambda})^2 \delta^{2\lambda - 1}\varphi^{(2)}_{n-1} \right\}.
\end{equation}
Equalities hold in (\ref{ce1})--(\ref{ce3}) only in the case $\theta =0$.
Further
\begin{eqnarray}
\label{e4}
{\varphi}_{n}^{(2)} & > & \delta ^{2\lambda -1}{\kappa }_{n}^{4}{\varphi}%
_{n-1}^{(2)}; \\
\label{e5}
{\varphi}_{n}^{(1)} & < & {\kappa }_{n}{\varphi}_{n-1}^{(1)};
\end{eqnarray}
\begin{equation}
\label{e6}
{\varphi}_{n}^{(1)}>{\kappa }_{n}{\varphi}_{n-1}^{(1)}+(\theta
+1)(1-\delta ^{-\lambda })\delta ^{2\lambda -1}{\kappa }_{n}^{3}\hat{\varphi}%
_{n-1}^{(2)};
\end{equation}
\begin{eqnarray}
\label{e1}
\tilde{\varphi}_{n}^{(2)} & > &\delta ^{-1}{\nu }_{n}^{4}\tilde{\varphi}%
_{n-1}^{(2)}; \\
\label{e2}
\tilde{\varphi}_{n}^{(1)} & < & {\nu }_{n}\tilde{\varphi}_{n-1}^{(1)};
\end{eqnarray}
\begin{equation}
\tilde{\varphi}_{n}^{(1)}>{\nu }_{n}\tilde{\varphi}_{n-1}^{(1)}+(\theta
+1)(1-\delta ^{-\lambda })\delta ^{-1}{\nu }_{n}^{3}\delta ^{-(n-1)\lambda }%
\tilde{\varphi}_{n-1}^{(2)}.  \label{e3}
\end{equation}
\end{Lm}

 {\bf Proof.} First we prove (\ref{ce1}). Consider
\begin{equation}
\label{ce4}
S(t,z) \ \stackrel{\rm def}{=} \ \exp(t\Delta_\theta )\left[g_{n-1} (z%
\delta^{-1-\lambda})\right]^\delta , \ \ t\in[0, \bar{t}], \ \ \bar{t}%
\stackrel{\rm def}{=} \delta^\lambda -1 , \ \ n\in{I\!\! N}.
\end{equation}
Taking into account Corollary \ref{XX1co}, (\ref{XX2}), (\ref{Z2}), and Lemma
\ref{Y1lm} one concludes that $S$ belongs to ${\cal L}^{+}$ and satisfies the
equation
\begin{eqnarray}
\label{ce5}
\frac{\partial S}{\partial t}  = \Delta_\theta S, \ \
S(0,z)  =  \left[g_{n-1} (z%
\delta^{-1-\lambda})\right]^\delta , \ \ S(\bar{t} ,z)
= g_n (z) .
\end{eqnarray}
We set
\[
%\label{ce7}
S_k (t) \ \stackrel{\rm def }{=} (D^k_z S)(t,0) , \ \ \ k\in{I\!\! N}_0 ,
\]
and obtain from (\ref{ce5}) and (\ref{a1})
\[
\frac{\partial S_0 (t)}{\partial t} = \theta S_1 (t) .
\]
Since $S\in {\cal L}^+$, $S_1 (t) > 0$ and for all $t\in [0,\bar{t}]$,
\[
S_0 (\bar{t}) > S_0 (0) \ \ {\rm for} \ \theta >0, \ \ S_0 (\bar{t}) = S_0 (0)
\ \ {\rm for} \ \theta = 0.
\]
This estimate and the boundary conditions (\ref{ce5}) gives (\ref{ce1}). Now we
set
\begin{equation}
\label{Z10}
p_{n}(z)=\exp (-{\varphi}_{n}^{(1)}z)g_{n}(z),
\end{equation}
insert $g_{n-1}(z) = \exp ({\varphi}%
_{n-1}^{(1)}z)p_{n-1}(z)$ into (\ref{Z2}), and use (\ref{ident}). Then
\begin{equation}
\label{M}
g_{n}(z)=\nu _{n}^{\theta }\exp (\kappa _{n}{\varphi}%
_{n-1}^{(1)}z)\exp (t_{n}\kappa _{n}^{-2}{\Delta }_{\theta
})\left[p_{n-1}(z\delta ^{\lambda -1}\kappa _{n}^{2})\right]^{\delta },
\end{equation}
where $t_{n}=(1-\delta ^{-\lambda })\kappa _{n}.$ For $t\in \lbrack 0,t_{n}]$,
we define
\begin{equation}
\label{R}
\exp R(t,z)=\exp (t{\Delta }_{\theta })\left[p_{n-1}(z\delta ^{\lambda
-1})\right]^{\delta }.
\end{equation}
According to Proposition \ref{Shift}, the above function can be written in the
form
\[
\exp R(t,z)=\exp (\hat{u}z)\hat{p}(z),
\]
where $\hat{u}<0$ and $\hat{p}$ belongs to  ${\cal L}^{+}$. Thus Proposition
\ref{SRule} yields for $k\geq 2$
\begin{equation}
\label{sr}
(-1)^{k-1}R_{k}(t) = (-1)^{k-1} (D^k_z \log \hat{p}(0)>0.
\end{equation}
Besides, we have
\begin{equation}
\label{ZZ1}
R(0,z)=\delta\log C_{n-1} -\delta ^{\lambda }{\varphi}_{n-1}^{(1)}z+\delta {\varphi}%
_{n-1}(z\delta ^{\lambda -1}),
\end{equation}
which gives
\begin{equation}
R_0 (0) = \delta\log C_{n-1},  \quad R_{1}(0)=0, \quad R_{2}(0)=\delta
^{2\lambda
-1}{\varphi}_{n-1}^{(2)}.
\label{CC1}
\end{equation}
Comparing (\ref{R}) and (\ref{M}), one obtains
\begin{eqnarray}
\label{Z12}
R(t_{n},z\kappa _{n}^{2})  =  {\varphi}_{n}(z)-\kappa _{n}{\varphi}%
_{n-1}^{(1)}z + \log C_{n} - \theta \log \nu_n ,
\end{eqnarray}
which yields
\begin{eqnarray}
\label{C2}
 &  & R_0 (t_n )  =  \log C_n - \theta \log \nu_n , \quad
R_{1}(t_{n})  =  \kappa _{n}^{-2}({\varphi}_{n}^{(1)}-\kappa _{n}{\varphi%
}_{n-1}^{(1)}), \\
  &  & R_{2}(t_{n} )  =  \kappa _{n}^{-4}{\varphi}_{n}^{(2)}. \nonumber
\end{eqnarray}
For $R(t,z)$, we obtain from (\ref{R}) an equation of the type of (\ref{ee1}),
(\ref{ce5})
\[
\frac{\partial R(t,z)}{\partial t}=\theta
(D_{z}R)(t,z)+z[(D_{z}^{2}R)(t,z)+(D_{z}R)^{2}(t,z)],
\]
with the initial condition given by (\ref{ZZ1}). It yields in turn
\begin{eqnarray}
\label{C13}
\frac{\partial R_0 (t) }{\partial t} & = & \theta R_1 (t), \\
\label{Z13}
\frac{\partial R_{1}(t)}{\partial t} & = & (\theta +1)R_{2}(t)+R_{1}^{2}(t), \\
\label{Z14}
\frac{\partial R_{2}(t)}{\partial t} & = & (\theta +2)R_{3}(t)+4R_{1}(t)R_{2}(t).
\end{eqnarray}
By the sign rule (\ref{sr}), $R_2 (t) <0$, thus for every $t_* $ such that $R_1
(t_* ) = 0$, one has from (\ref{Z13})
\[
\frac{\partial R_{1}}{\partial t} (t_* )< 0.
\]
Clearly, such $t_* $ is at most one. Since $R_{1}(0)=0$ , one has $t_* = 0$ and
\begin{equation}
\label{Z15}
R_{1}(t)<0,\ \ \forall t\in (0,t_{n}],
\end{equation}
which yields in (\ref{C13})
\[
R_0 (t_n ) > R_0 (0) \ \ {\rm for} \ \theta>0, \ \ \ R_0 (t_n ) = R_0 (0) \ \
{\rm for} \ \theta=0 ,
\]
and $ R_1 (t_n ) <0 $, implying (\ref{ce2}) and (\ref{e5}) if the conditions
(\ref{ZZ1}) --  (\ref{C2}) are taken into account. Applying again (\ref{sr}) and
(\ref{Z15}) in (\ref{Z14}) we get
\begin{equation}
\frac{\partial R_{2}(t)}{\partial t}>0,\ \ \forall t\in (0,t_{n}],
\label{R2t}
\end{equation}
which yields in (\ref{Z13})
\begin{equation}
\label{CC10}
R_1 (t) > t(\theta+1)R_2 (0)
\end{equation}
and
\begin{equation}
\label{Es1}
 R_{2}(0)<R_{2}(t_{n}).
\end{equation}
The latter gives (\ref{e4}).  Taking in (\ref{CC10}) $t=t_n $ one obtains (\ref{e6}).
 To obtain (\ref{ce3}) one observes that (\ref{CC10}) and (\ref{Z15}) yield in (\ref{C13})
for $\theta >0$
\[
R_0 (t_n ) - R_0 (0) > \frac{1}{2}t_n^2 \theta (\theta +1) R_2 (0).
\]
For $\theta =0$, we have already obtained $R_0 (t_n ) = R_0 (0)$.  Finally,
(\ref{e1})-(\ref{e3}) follow directly from (\ref{e4})-(\ref{e6}) and (\ref{Z8}
).\qquad \hfil\vrule height6pt width5pt depth-1pt
\newline By the first condition in (\ref{ee8}), there exists $\sigma \in [\delta
^{(2\lambda -1)/4},1)$ such that
\begin{equation}
\frac{m_{2}(g)}{[\alpha +m_{1} (g)]^{2}}=\frac{\delta ^{1/2}}{\theta +1}\frac{%
1-\sigma }{\delta ^{\lambda }-\sigma }.  \label{con2}
\end{equation}
For such $\sigma $, we set
\begin{equation}
\Phi ^{(1)} \ \stackrel{\rm def}{=} \ \frac{1-\sigma \delta ^{-\lambda }}{1-\delta ^{-\lambda }},
\label{Fi1}
\end{equation}
\begin{equation}
\label{Fi2}
\Phi ^{(2)}\ \stackrel{\rm def}{=} \ -\Phi ^{(1)}
\frac{\delta ^{1-\lambda }}{\theta +1}\frac{\sigma^{2}(1-\sigma )}
{\delta ^{\lambda }-1}.
\end{equation}
\begin{Lm}
\label{Triple}The following triple ${\cal I}_{n}=\left(
i_{n}^{1};i_{n}^{2};i_{n}^{3}\right) $\ of statements:
\begin{eqnarray*}
i_{n}^{1} & = & \left\{ \exists \beta _{n}^{+}\in J(g):%
{\varphi}_{n}^{\left( 1\right) }=\Phi ^{\left( 1\right) },\beta =\beta
_{n}^{+};\quad {\varphi}_{n}^{\left( 1\right) }<\Phi ^{\left( 1\right)
}, \ \beta <\beta _{n}^{+}\right\} , \\
i_{n}^{2} & = & \left\{ \exists \beta _{n}^{-}\in J(g):%
{\varphi}_{n}^{\left( 1\right) }=1,\beta =\beta _{n}^{-};\quad {\varphi}%
_{n}^{\left( 1\right) }<1, \ \beta <\beta _{n}^{-}\right\} , \\
i_{n}^{3} & = & \left\{ \forall \beta \leq \beta _{n}^{+}:{\varphi}%
_{n}^{\left( 2\right) }\geq \Phi ^{\left( 2\right) }\right\},
\end{eqnarray*}
holds true for all $n\in {I\!\!N}_{0}.$
\end{Lm}

{\bf Proof.} For $n=0$, we have ${\varphi}_{0}^{(1)}=\beta (\alpha
+m_{1} (g))$, ${\varphi}_{0}^{(2)}=-\beta ^{2}m_{2} (g)$. Thus we set
\begin{equation}
\label{beta}
\beta _{0}^{-}=\frac{1}{\alpha +m_{1} (g)},\quad \beta _{0}^{+}=\frac{\Phi
^{\left( 1\right) }}{\alpha +m_{1}(g)}>\beta _{0}^{-} .
\end{equation}
First let us prove that $\beta _{0}^{+}\in J(g)$. If $\alpha =0,$ $\beta
_{0}^{+}$ needs only to be finite, which obviously holds. For $\alpha >0$, the definitions (\ref{beta}) and (\ref
{Fi1}) yield for $\beta =\beta _{0}^{+}$
\[
{\varphi}_{0}^{(1)}=\Phi ^{\left( 1\right)
}=\frac{\delta ^{1/2}}{\delta ^{1/2}-(\theta +1)m_{2}(g)/[\alpha +m_{1}(g)]^{2}},
\]
thus
\[
{\varphi}_{0}^{(1)}=\frac{\delta ^{1/2}}{\delta ^{1/2}-(\theta +1)(\beta
_{0}^{+})^{2}m_{2}(g)/({\varphi}_{0}^{(1)})^{2}}.
\]
This equation can be solved with respect to ${\varphi}_{0}^{(1)}$%
\[
{\varphi}_{0}^{(1)}=\frac{1}{2}\{1+[1+4\delta ^{-1/2}(\theta +1)(\beta
_{0}^{+})^{2}m_{2}(g)]^{1/2}\}.
\]
Hence making use of the second condition in (\ref{ee8}) one gets
\begin{eqnarray*}
\beta _{0}^{+}(\alpha +m_{1}(g)) & = & {\varphi}_{0}^{(1)}<\frac{1}{2}\{%
1+1+[4\delta ^{-1/2}(\theta +1)(\beta _{0}^{+})^{2}m_{2}(g)]^{1/2}\} \\
& = & 1+[\delta ^{-1/2}(\theta +1)(\beta _{0}^{+})^{2}m_{2} (g)]^{1/2}\leq 1+\beta
_{0}^{+}m_{1} (g).
\end{eqnarray*}
Therefore, $\beta _{0}^{+}\in J(g)$ and $ i_{0}^{1},$ $i_{0}^{2}$ are true. To
prove $i_{0}^{3}$ we to apply ( \ref{con2}). Indeed, for $\beta=\beta_{0}^{+}$,
\begin{eqnarray*}
{\varphi}_{0}^{(2)} & = & -(\beta _{0}^{+})^{2}m_{2}(g) =-(\Phi ^{\left( 1\right)
})^{2}\frac{m_{2}(g)}{[\alpha +m_{1}(g)]^{2}}=-(\Phi ^{\left( 1\right) })^{2}\frac{%
\delta ^{1/2}}{\theta +1}\frac{1-\sigma }{\delta ^{\lambda }-\sigma } \\
& = &  (\Phi ^{\left( 1\right) })^{2}(\Phi ^{\left( 1\right) })^{-2}\Phi
^{(2)}\delta ^{\lambda -1/2}\sigma ^{2}\geq \Phi ^{(2)}.
\end{eqnarray*}
This proves ${\cal I}_{0}$. Note that the estimate (\ref{Ineq}) with $n=0$ holds
for $\beta \in (0, \beta^+_{0}]$. To prove the implication ${\cal
I}_{n-1}\Rightarrow {\cal I}_{n}$, we remark that, for $\beta
=\beta_{n-1}^{+}$, $i_{n-1}^{1}$ yields
${\varphi}_{n-1}^{(1)}=\Phi ^{(1)}$ and ${\kappa }_{n}=\sigma ^{-1}$ ( see (%
\ref{kap})). Now for $\varphi_n^{(1)}$, we have the following possibilities: (a)
the estimate (\ref{Ineq}) holds for $\beta = \beta_{n-1}^+$; (b) this estimate
does not hold for such $\beta$. In the case (a) one may apply Lemma \ref{Z2lm}.
Then by means of $i_{n-1}^{3}$, (\ref{e6}), (\ref{Fi1}), and (% \ref{Fi2}), we
obtain
\[
{\varphi}_{n}^{(1)}>\sigma ^{-1}\Phi ^{(1)}+(\theta +1)(1-\delta
^{-\lambda })\sigma ^{-3}\delta ^{2\lambda -1}\Phi ^{(2)}=\Phi ^{(1)}.
\]
In the case (b) we simply have
\[
{\varphi}_{n}^{(1)}\geq \frac{1}{1-\delta ^{-\lambda }}\geq \frac{%
1-\sigma \delta ^{-\lambda }}{1-\delta ^{-\lambda }}=\Phi ^{(1)}.
\]
For $\beta =\beta _{n-1}^{-}$, we have ${\varphi}_{n-1}^{(1)}=1$ and ${%
\kappa }_{n}=1$. Therefore , ${\varphi}_{n}^{(1)}<1$ for $\beta \leq
\beta _{n-1}^{-}$, as follows from $i_{n-1}^{2}$ and (\ref{e5}). By Corollary
\ref{DiffB}, ${\varphi}_{n}^{(1)}$ is a continuous function of $%
\beta $, thus there exists at least one value of $\beta =\tilde{\beta}%
_{n}^{+}\in (\beta _{n-1}^{-},\beta _{n-1}^{+})$ such that ${\varphi}%
_{n}^{(1)}=\Phi ^{(1)}$. The smallest such one is set to be ${\beta }%
_{n}^{+} $. The existence of ${\beta }_{n}^{-}\in (\beta _{n-1}^{-},{\beta }%
_{n}^{+})$ can be established in the same way. For $\beta \leq {\beta }%
_{n}^{+}$, we have $\beta \leq \beta _{n-1}^{+}$ and then ${\varphi}%
_{n-1}^{(1)}\leq \Phi ^{(1)}$ due to $i_{n-1}^{1}$. This yields ${\kappa }%
_{n}\leq \sigma ^{-1}$, then we get from (\ref{e4})
\[
{\varphi}_{n}^{(2)}>\sigma ^{-4}\delta ^{2\lambda -1}{\varphi}%
_{n-1}^{(2)}\geq {\varphi}_{n-1}^{(2)}\geq \Phi ^{(2)},
\]
where the following estimates were used: ${\varphi}_{n}^{(2)}<0$, $%
\forall n\in {I\!\! N}_{0}$ ; $\sigma ^{-1}\leq \delta ^{(1-2\lambda )/4}$.
\qquad \hfil\vrule height6pt width5pt depth-1pt
\begin{Co}
\label{Ineqco}
The inequality (\ref{Ineq}) holds for $\beta \leq \beta_n^+$, thus $J_n $ is
nonempty.
\end{Co}
\begin{Lm}
\label{Bstar}
There exists $\beta _{\ast }\in J(g)$
such that, for $\beta =\beta _{\ast }$,
\begin{equation}
1<{\varphi}_{n}^{(1)}<\Phi ^{(1)},\ \forall n\in {I\!\!N}_{0}.
\label{e7}
\end{equation}
For $\beta <\beta _{\ast }$, the above upper estimate also holds and, moreover,
there exists $K=K(\beta )>0$, such that
\begin{equation}
{\varphi}_{n}^{(1)}<K\delta ^{-\lambda n},\ \forall n\in {I\!\!N}_{0}.
\label{e8}
\end{equation}
\end{Lm}

{\bf Proof. } Consider the set $\Delta _{n}\ \stackrel{\rm def}{=}%
\{\beta \in (0,{\beta }%
_{n}^{+})\mid 1<{\varphi}_{n}^{(1)}<\Phi ^{(1)}\}$.
Just above we have shown that
$\Delta _{n}\subseteq ({\beta }_{n}^{-},{\beta }_{n}^{+})$, $\Delta
_{n}$ is nonempty and open. Let us prove that $\Delta _{n}\subseteq \Delta
_{n-1}$. Suppose there exists some $\beta \in \Delta _{n}$ which does not
belong to $\Delta _{n-1}$. For this $\beta $, either ${\varphi}%
_{n-1}^{(1)}\leq 1$ or ${\varphi}_{n-1}^{(1)}\geq \Phi ^{(1)}$. Hence
either ${\varphi}_{n}^{(1)}<1$ or ${\varphi}_{n}^{(1)}>\Phi ^{(1)}$ (it can be
proved as above). This runs in counter with the supposition $\beta
\in \Delta _{n}$, hence $\Delta _{n}\subseteq \Delta _{n-1}.$ Now let
$D_{n}$ be the closure of $\Delta _{n}$, then
\begin{equation}
D_{n}=\{\beta \in \lbrack {\beta }_{n}^{-},{\beta }_{n}^{+}]\mid 1\leq
\varphi_{n}^{(1)}\leq \Phi ^{(1)}\}.  \label{Dn}
\end{equation}
$D_{n}$ is nonempty and $D_{n}\subseteq D_{n-1}\subseteq ...\subseteq
D_{0}\subseteq J(g)$. Let $D_{\ast }=\bigcap_{n}D_{n}$%
, then $D_{\ast }$ is also nonempty and closed, and $D_{\ast }\subset J(g).$
Now let us show that, for every $\beta \in D_{\ast },$ the
estimates (\ref{e7}) hold. Indeed, directly from the definition of $D_{\ast
} $ one has
\[
1\leq {\varphi}_{n}^{(1)}\leq \Phi ^{(1)},\ \forall n\in {I\!\!N}_{0}.
\]
Suppose ${\varphi}_{n}^{(1)}=1$ for some $n\in {I\!\!N}_{0}$ and $\beta
\in D_{\ast }$, then ${\varphi}_{m}^{(1)}<1$ for all $m>n$ (see (\ref{e5}%
)). The latter means that this $\beta $ does not belong to all $D_{m}$ with
$m>n$. This contradicts the supposition $\beta \in D_{\ast }.$ The case
${\varphi}_{n}^{(1)}=\Phi ^{(1)}$ can be excluded similarly. Set $\beta
_{\ast }=\min D_{\ast }$. We have just proved that, for $\beta
=\beta _{\ast }$, (\ref{e7}) holds, thus it remains to prove the second
part of Lemma. To this end we take $\beta <\beta _{*}.$ If
${\varphi}_{n}^{(1)}>1$ for all $ n\in {I\!\! N}_{0},$ then either (\ref{e7})
holds or there exists such $n_{0}$ that ${\varphi}_{n_{0}}^{(1)}\geq \Phi
^{(1)}.$ This means either $\beta \in D_{*}$ or $\beta >\inf {\beta }_{n}^{+}.$
Both these cases contradict the definition of $\beta _{*}$. Hence there exists $
n_{0}$ such that ${\varphi}_{n_{0}-1}^{(1)}\leq 1,$ then ${\varphi}%
_{n}^{(1)}<1$ for all $n\geq n_{0}.$ In what follows, the definition (\ref
{kap}) and the estimate (\ref{e5}) imply for the sequences $\{{\varphi}%
_{n}^{(1)}, \ n\geq n_{0}\}$ and $\{\kappa _{n}, \ n\geq n_{0}\} $ to be
strictly decreasing. Then for all $n>n_{0}$, one has (see (\ref{e5}))
\[
{\varphi}_{n}^{(1)}<\kappa _{n}{\varphi}_{n-1}^{(1)}<...<\kappa
_{n}\kappa _{n-1}...\kappa _{n_{0}+1}{\varphi}_{n_{0}}^{(1)}<(\kappa
_{n_{0}+1})^{n-n_{0}}.
\]
Since $\kappa _{n_{0}+1}<1,$ one has $\sum_{n=0}^{\infty }%
{\varphi}_{n}^{(1)}<\infty
.$ Thus there exists $0<K_{0}<\infty $ such that
\begin{equation}
\prod_{n=1}^{\infty }\nu _{n} \stackrel{\rm def}{=}K_{0}.  \label{C0}
\end{equation}
Finally, we apply (\ref{e5}) once again and obtain
\begin{equation}
\label{C1}
{\varphi}_{n}^{(1)}<\delta ^{-\lambda n}\nu _{n}\nu _{n-1}...\nu _{1}%
{\varphi}_{0}^{(1)}<\delta ^{-\lambda n}K_{0}{\varphi}%
_{0}^{(1)} \ \stackrel{\rm def}{=} \ K\delta ^{-\lambda n},\ \forall n\in {I\!\! N}_{0}.\qquad
\end{equation}
\hfil\vrule height6pt width5pt depth-1pt
\newline Now we state the lemmas the proof of our theorems directly follows from. The
first four lemmas describe the sequences $\{g_n \}$ defined by (\ref{Z2}) whose
elements have the form (\ref{Z7b}).
\begin{Lm}
\label{Zb11lm}
For every $\theta \geq 0$ and $g\in {\cal L} (\lambda)$, there exists $\beta_*
\in J(g)$ such that,
\vskip.1cm
\begin{tabular}{ll}
{\rm (i)} &  for
$\beta = \beta_* $, $\lim_{n\rightarrow \infty} \varphi_n^{(1)} =1 $ and
$\ \lim_{n\rightarrow \infty} \varphi_n^{(2)} =0 $; \\[.1cm]
{\rm (ii)} & for
$\beta < \beta_* $, $\lim_{n\rightarrow \infty} \varphi_n^{(1)} = 0 $ and
$\ \lim_{n\rightarrow \infty} \varphi_n^{(2)} =0 $.
\end{tabular}
\end{Lm}
\begin{Lm}
\label{Za11lm}
Let $\theta$, $g$ and $\beta_* $ be as above. Then there exists $C:(0,\beta_* ]
\rightarrow {I\!\! R}_+ $ such that the sequence $\{ C_n \ | \  n\in {I\!\! N}_0 , \
C_n = g_n (0), \ C_0 =  C(\beta ) \}$, converges to $C_{2,*}$ (resp. to
$C_{1,*}$) given by (\ref{Y2}) for $\beta = \beta_*$ (resp. $\beta <\beta_* $).
The sequence $\{ C_n \ | \ C_0 > C(\beta) \}$ is divergent, the sequence $\{ C_n
\ |
\
 C_0 < C(\beta) \}$ tends to zero.
\end{Lm}
\begin{Lm}
\label{Z11lm}
Let $\theta$, $g$, $\beta_* $ and $C(\beta)$ be as above. Then for $\beta =
\beta_* $, the sequence $\{ g_n \ | \ n\in {I\!\! N}_0 , \ g_0 (z) = C(\beta_*
)g(\beta z)\}$ converges in ${\cal A}_1 $ to $g_{2,*} (z) = C_{2,*}\exp(z)$
defined by (\ref{Y3a}).
\end{Lm}
\begin{Lm}
\label{Z1lm}
Let $\theta$, $g$, $\beta_* $ and $C(\beta)$ be as above. Then for every $\beta
< \beta_*
$
\vskip.1cm
\begin{tabular}{ll}
{\rm (i)}  & the sequence $\{ g_n \ | \ n\in {I\!\! N}_0 , \ g_0 (z) = C(\beta
)g (\beta z) \}$ converges in ${\cal A}_1$ \\ & to $g_{1, *}(z) \equiv 1$ ;
\\[.1cm] {\rm (ii)} & the sequence $\{ \tilde{g}_n \ | \ n\in {I\!\! N}_0 , \
\tilde{g}_0 (z)%
 = g (\beta z) \}$ defined by (\ref{Z7a}) \\
  &  converges in ${\cal A}_\varphi $ to $\tilde{g}_{*} (z) = \exp(\varphi z )$
     with certain $\varphi = \varphi (\beta)>0 $.
\end{tabular}
\end{Lm}
Directly from the definitions (\ref{23}), (\ref{24}), and (\ref{Z2}) one has the
following corollary of the above lemmas.
\begin{Lm}
\label{Z13lm}
For every $\theta \geq 0$ and $g\in {\cal L}(\lambda )$, there exist $\tau_* \in
I(g)$ and a function $C: [0.\tau_*] \rightarrow {I\!\! R}_+ $, such that:
\vskip.1cm
\begin{tabular}{ll}
{\rm (i)}   &  for $\tau < \tau_* $, the sequence
$\{ f_n \ | \ n\in {I\!\! N}_0 , \ f_0 (z)  = C(\tau ) g (z)  \}$ \\
& converges in ${\cal A}_{\beta_*^{-1} }$
  to $f_{1,*} (z) \equiv 1 $; \\[.1cm]
{\rm (ii)}& for $\tau = \tau_* $,
the sequence $\{ f_n \ | \ n\in {I\!\! N}_0 , \ f_0 (z) = C(\tau_* ) g (z)  \}$ \\
& defined
by (\ref{23}), (\ref{24}), and (\ref{X5}) converges in ${\cal A}_{\beta_*^{-1} }$, \\
& $\beta_* = \tau_* (\delta^{\lambda} -1)^{-1}$
 to  $f_{2,*} (z)
= \delta^{-\lambda \theta /(\delta - 1)}\exp(\beta_*^{-1} z)$.
\end{tabular}
\end{Lm}

\subsection{Proof of Theorems}

        {\bf Proof of Theorem \ref{ee0tm}. }As it has already been established, the function
$f_n (t, z)$ defined by (\ref{23}) gives the solution of the problem (\ref{ee5}) provided
all $f_{m}$, $m= 0,1, \dots n-1$ belong to the domain of the operators $T_t $, $t\in [0,1]$.
 The latter fact follows from Lemma \ref{Y1lm} and Proposition \ref{5Xpn}.
\hfil\vrule height6pt width5pt depth-1pt
\newline {\bf Proof of Theorem \ref{ee1tm}. }Proposition \ref{X4pn} and Lemma \ref{Z13lm}
yield that, for $\tau < \tau_*$, the sequence $\{f_n (t,z) \}$ converges in
${\cal A}_{\beta_*^{-1}}$ to $T_t (f_{1,*})$, which may be easily calculated to
be identically one. For $\tau = \tau_* $, one has the same convergence to $T_t%
(f_{2,*})$, which can be calculated explicitly by means of Proposition
\ref{Shift}.
\hfil\vrule height6pt width5pt depth-1pt
\newline The proof of Theorems \ref{2tm} and \ref{3tm}
follows directly from Theorems \ref{ee0tm} and \ref{ee1tm} respectively on the
base of the identity (\ref{17}).
\newline {\bf Proof of Theorem \ref{4tm} and Theorem \ref{5tm}. }By the continuity theorem
(see e.g. \cite{[Ibr]}, p.27), the convergence of the sequence of the transforms
(\ref{9y}) $\{F_{\mu_n }\}$ in a certain ${\cal A}^{(N)}_a $ implies the weak
convergence of the sequence $\{\mu_n
\}$. But for the measures defined by (\ref{11y}), the transforms
(\ref{9y}) are the isotropic functions $F_{\mu_n }\in {\cal A}^{(N)}_a $ and for
any of them, there exists an entire function $f_{\mu_n } \in {\cal A}_a $
obeying (\ref{11x}). Moreover, (\ref{11y}) implies
\[
f_{\mu_n } (z) = T(f_{\mu_{n-1}} ) (z)\left[T(f_{\mu_{n-1}} ) (0) \right]^{-1},
\ \
\ n\in {I\!\! N}.
\]
Now if one chooses the starting element $\mu_0 = \nu$ such that
\[
\int_{{I\!\! R}^N } \exp((x,y))\nu(dy) = g((x,x)),
\]
where $g $ is the starting element of $\{{f}_n \}$ described by Lemma
\ref{Z13lm}, then the validity of Theorem \ref{4tm} follows from this Lemma. The
assertion regarding the variance in claim (ii) may be checked directly.
Similarly, the transforms $F_{\tilde{\mu}_n }$ (\ref{9y}) of $\tilde{\mu}_n  $,
defined by (\ref{15y}), and the elements of the sequence $\{
\tilde{g}_n \}$, defined by (\ref{Z7a}) and described by claim (ii) of Lemma
\ref{Z1lm}, obey the relation
\[
F_{\tilde{\mu}_n }(\sqrt{\beta}x ) = \tilde{g}_n ((x,x)), \ \  n\in {I\!\! N},
\ \ \ F_{\tilde{\mu}_0 }(x ) = g((x,x)) .
\]
Then the validity of Theorem \ref{5tm} follows directly from claim (ii) of Lemma \ref{Z1lm}.
\hfil\vrule height6pt width5pt depth-1pt

\subsection{Proof of Lemmas}

{\bf Proof of Lemma \ref{Zb11lm}.} Consider the case $\beta =\beta _{\ast },$
where (\ref{e7}) holds and ${\varphi}_{0}^{(2)}\geq \Phi ^{(2)}$ by statement $
i^3_0 $ of Lemma \ref{Triple} . First we prove that
${\varphi}_{n}^{(2)}\rightarrow 0\ $. Here we have such two possibilities:

${\rm (a)} \quad \sigma >\delta ^{(2\lambda -1)/4}.\quad $From (\ref{Z9}) and
(\ref {e7}) we obtain ${\kappa }_{n}<\sigma ^{-1}$. Thus
\begin{equation}
\delta ^{2\lambda -1}{\kappa }_{n}^{4}<\delta ^{2\lambda -1}\sigma ^{-4}<1
\label{n1}
\end{equation}
Applying this estimate in (\ref{e4}) one gets
\[
\mid {\varphi}_{n}^{(2)}\mid <\delta ^{2\lambda -1}\sigma ^{-4}|\hat{%
\varphi}_{n-1}^{(2)}|<...<(\delta ^{2\lambda -1}\sigma ^{-4})^{n}\mid \hat{%
\varphi}_{0}^{(2)}\mid \leq (\delta ^{2\lambda -1}\sigma ^{-4})^{n}\mid \Phi
^{(2)}\mid .
\]
In view of (\ref{n1}), this gives
\[
{\varphi}_{n}^{(2)}\rightarrow 0,\ \ n\rightarrow +\infty .
\]

${\rm (b)}\quad \sigma =\delta ^{(2\lambda -1)/4}.\quad $In this case we have only
\begin{equation}
\delta ^{2\lambda -1}{\kappa }_{n}^{4}<\delta ^{2\lambda -1}\sigma ^{-4}=1.
\label{n2}
\end{equation}
Making use of (\ref{e4}) one obtains
\[
0>{\varphi}_{n}^{(2)}>\delta ^{2\lambda -1}{\kappa }_{n}^{4}{\varphi}%
_{n-1}^{(2)}>{\varphi}_{n-1}^{(2)}>...>{\varphi}_{0}^{(2)}\geq \Phi
^{(2)}
\]
Hence $\{{\varphi}_{n}^{(2)}\}$ is strictly increasing and bounded. Then it is
convergent and its limit, say ${\varphi}^{(2)},$ obeys the condition $
{\varphi}^{(2)}>\Phi ^{(2)}$. Assume now that ${\varphi}^{(2)}\neq 0.$ Combining
(\ref{e4}) and (\ref{n2}) one obtains (recall that ${\varphi}_{n}^{(2)}<0$)
\[
\frac{{\varphi}_{n}^{(2)}}{{\varphi}_{n-1}^{(2)}}<\delta ^{2\lambda
-1}{\kappa }_{n}^{4}<1,
\]
which means ${\kappa }_{n}\rightarrow \delta ^{(1-2\lambda )/4}.$ The
latter as well as the definitions of ${\kappa }_{n}$ and $\Phi ^{(1)}$ immediately
yield
\[
{\varphi}_{n}^{(1)}\rightarrow \Phi ^{(1)}.
\]
Passing to the limit $n\rightarrow +\infty $ in (\ref{e6}) one obtains
\[
\Phi ^{(1)}\geq \delta ^{(1-2\lambda )/4}\Phi ^{(1)}+(\theta +1)(1-\delta
^{-\lambda })\delta ^{2\lambda -1}\delta ^{3(1-2\lambda )/4}{\varphi}%
^{(2)},
\]
which yields in turn
\[
{\varphi}^{(2)}\leq -\Phi ^{(1)}\frac{(1-\sigma )\sigma
^{2}\delta ^{1-\lambda }}{(\theta +1)(\delta ^{\lambda }-1)}=\Phi ^{(2)}.
\]
The latter gives the following contradictory inequalities
\[
\Phi ^{(2)}<{\varphi}^{(2)}\leq \Phi ^{(2)}.
\]
Thus ${\varphi}^{(2)}=0$. To show that ${\varphi}_{n}^{(1)}\rightarrow 1$, we
set
\begin{equation}
\label{be1}
b_{n}=(\theta +1)(1-\delta ^{-\lambda })\delta ^{2\lambda -1}{\kappa }%
_{n}^{3}{\varphi}_{n-1}^{(2)}.
\end{equation}
Since $\{{\kappa }_{n}\}$ is bounded and ${\varphi}_{n}^{(2)}\rightarrow 0,$ one
has $b_{n}\rightarrow 0.$ By the estimate ( \ref{e7}), $\{{\varphi}_{n}^{(1)}\}$
is bounded. Then it contains a subsequence $\{{\varphi}_{n_{i}}^{(1)}\}$
convergent to a certain $a\in [1,
\Phi^{(1)}]$. From (\ref{e5}) and (\ref{e6}) one has
\[
0>{\varphi}_{n}^{(1)}-{\kappa }_{n}{\varphi}_{n-1}^{(1)}>b_{n},
\]
which yields
\[
\lim_{i\rightarrow \infty }({\varphi}_{n_{i}}^{(1)}-{\kappa }_{n_{i}}%
{\varphi}_{n_{i}-1}^{(1)})=0.
\]
The latter can be rewritten as (see \ref{Z9})
\[
a-\frac{\delta ^{-\lambda }a}{1-(1-\delta ^{-\lambda })a}=0.
\]
Since $\lambda >0,$ the above equation has only one solution on $\left[ 1,\Phi
^{(1)}\right] $, it is $a=1$. In what follows, the bounded sequence
$\{{\varphi}_{n}^{(1)}\}$ has only one accumulation point, hence it converges to
$a=1$ itself. In the case $\beta <\beta_*$ the estimate (\ref{e8}) yields
$\varphi_n^{(1)} \rightarrow 0$. Then $\kappa_n $ given by (\ref{Z9}) tends to
$\delta^{-\lambda}$ which immediately gives in (\ref{e4})
$\varphi_n^{(2)}\rightarrow 0$.
\qquad \hfil\vrule height6pt width5pt depth-1pt

{\bf Proof of Lemma \ref{Za11lm}.} From the definitions (\ref{Z2}) and
(\ref{Z7b}) one obtains
\begin{eqnarray*}
C_n & = &  C_{n-1}^\delta \Psi_{n-1} (\beta) , \\
\Psi_n (\beta) & \stackrel{\rm def}{=} & \left\{ \exp\left((\delta^\lambda -1)%
\Delta_\theta \right) \exp\left( \delta \varphi_n (\delta^{-1-\lambda} \cdot)\right)%
\right\} (0).
\end{eqnarray*}
For $\theta =0$, $\Psi_n (\beta) = 1$ (see Remark \ref{n2rk}) and the situation
with $C_n $ is obvious. Consider the case $\theta >0$. Then
\begin{equation}
\label{L1}
C_n = C_0^{\delta^n }\Xi_n (\beta),  \ \  \Xi_n (\beta) \ \stackrel{\rm def}{=}
\ \Psi_{n-1}(\beta) \Psi_{n-2}^\delta (\beta) \dots \Psi_0^{\delta^{n-1}}(\beta).
\end{equation}
Now we put $C_0 = \zeta >0$, then $C_n = C_n (\zeta , \beta)$. By the above representation,
for every fixed $\beta>0$, $C_n $ is a monotone convex differentiable function
of $\zeta$ and
\begin{equation}
\label{L2}
C_n (\zeta , \beta )= \zeta^{\delta^n}\Xi_n (\beta) , \ \
\frac{\partial C_n }{\partial \zeta } = \delta^{n}\zeta^{-1}C_n .
\end{equation}
By Lemma \ref{Bstar},  $\varphi_n^{(1)} <
\Phi^{(1)}$ for all $n\in {I\!\! N}_0 $ and $\beta \in (0, \beta_* ]$. This gives in (\ref{Z9})
$\kappa_n < \sigma^{-1}\leq \delta^{(1-2\lambda)/4}$ for such $\beta$ and $n$.
We set
\begin{eqnarray}
\label{L3}
\zeta^- & = & \delta^{-\theta (1+2\lambda) /4(\delta -1)}, \\
\label{L4}
\Upsilon & = & [\zeta^- , 1] \subset {I\!\! R}_+ .
\end{eqnarray}
For a fixed $\beta \in (0, \beta_* ]$, let us prove that the following inductive
statements hold true for all $n\in {I\!\! N}_0 $
\begin{eqnarray}
\label{L5}
i_n^+ & = & \{ \exists \zeta_n^+ \in \Upsilon : C_n (\zeta_n^+ , \beta) = 1  \}, \\
i_n^- & = & \{ \exists \zeta_n^- \in \Upsilon : C_n (\zeta_n^- , \beta) = \zeta^-  \}.
\nonumber
\end{eqnarray}
Since $C_n $ is a monotone convex function of $\zeta$ (\ref{L2}), such
$\zeta_n^{\pm}$ are unique. We set $\zeta_0^+ = 1$, $\zeta_0^- = \zeta^- $. Then
$C_0 = \zeta $ obeys the above conditions, thus $i_0^{\pm}$ are true. Now
suppose that $i_{n-1}^{\pm}$ are true. Then (\ref{ce1}) and (\ref{ce2}) yield
for $\theta >0$
\[
C_n (\zeta_{n-1}^+ , \beta ) > 1, \ \ \ C_n (\zeta_{n-1}^- , \beta ) <
\zeta^-
.
\]
Taking into account that $C_n $ depends on $\zeta$ as given by (\ref{L2}) one concludes
that there exist $\zeta_n^{\pm}$ such that
\begin{equation}
\label{L6}
\zeta_{n-1}^- < \zeta_n^- < \zeta_n^+ < \zeta_{n-1}^+ ,
\end{equation}
and the statements $i_n^{\pm}$ hold true. Set
\begin{equation}
\label{L7}
\Upsilon_n = [\zeta_n^- , \zeta_n^+ ] .
\end{equation}
Then
\[
\Upsilon_n \subset \Upsilon_{n-1} \subset \dots \subset \Upsilon ,
\]
and there exists $\tilde{\zeta}_n \in \Upsilon_n $ such that
\[
\zeta_n^+ - \zeta_n^- = (1-\zeta^- )\left[\frac{\partial C_n (\tilde{\zeta}_n , \beta )}%
{\partial \zeta }\right]^{-1}.
\]
Let
\[
\Upsilon^* = \bigcap_{n\in {I\!\! N}_0 } \Upsilon_n ,
\]
which is closed and nonempty. For $\zeta \in \Upsilon^*$, all $C_n $ belong to
$\Upsilon$. Hence the sequence $\{C_n \}$ is separated from zero for such
$\zeta$. This yields that the derivative given by (\ref{L2}) tends to $+\infty$
when $n\rightarrow \infty$. Taking into account all these facts one concludes
\begin{equation}
\label{L8}
\Upsilon^* = \{ \zeta^* \}, \ \ \zeta^* \in \Upsilon
\end{equation}
and, for all $n\in {I\!\! N}_0 $,
\begin{equation}
\label{L9}
C_n (\zeta^* , \beta ) \in \Upsilon .
\end{equation}
It should be pointed out that $\zeta^* = \zeta^* (\beta )$. Choose $\zeta
= \zeta^* $. Then by (\ref{L9}), the sequence $\{C_n \}$ is bounded, hence it
contains a convergent subsequence. For $\beta = \beta_*$, by means of
(\ref{ce2}) and (\ref{ce3}) one may show that such a subsequence converges to
$C_{2,*} = \delta^{-\lambda \theta /(\delta-1)}$. As in the case of
$\{\varphi_n^{(1)} \}$ considered above, this fact implies the convergence of
the whole sequence to this limit. For $\beta <\beta_*$, one employs (\ref{ce1})
and (\ref{ce2}) and shows similarly the convergence of $\{C_n \}$ to $C_{1,*} =
1$. Thus we choose the function $C(\beta)$ to be $C(\beta) = \zeta^* (\beta )$.
\qquad \hfil\vrule height6pt width5pt depth-1pt

{\bf Proof of Lemma \ref{Z11lm}.} It follows from Lemmas \ref{Conver},
\ref{Zb11lm}, and \ref{Za11lm}. \qquad \hfil\vrule height6pt width5pt depth-1pt

{\bf Proof of Lemma \ref{Z1lm}.} Claim (i) follows from the lemmas just
mentioned. To prove claim (ii) we fix $\beta <\beta _{\ast }$ and show the
convergence of $\left\{
\tilde{\varphi}_n^{(2)}\right\} $ to zero. Indeed, (\ref{e1}) and (\ref{C0})
imply
\begin{equation}
\left| \tilde{\varphi}_{n}^{(2)}\right| <\delta ^{-n}(\nu _{n}\nu
_{n-1}...\nu _{1})^{4}\left| \tilde{\varphi}_{0}^{(2)}\right| <\delta
^{-}{}^{n}K_{0}^{4}\left| \tilde{\varphi}_{0}^{(2)}\right| .  \label{efi2}
\end{equation}
Thus to complete the proof we have only to show that $\{\tilde{\varphi}%
_{n}^{(1)}\}$ is a Cauchy sequence. To this end for $n\in {I\!\!N}$ and $%
p\in {I\!\!N}$, we set
\begin{equation}
a_{n,p}={\nu }_{n+p}\nu _{n+p-1}...\nu _{n+1};
\end{equation}
\begin{equation}
b_{n,p}=(\theta +1)(1-\delta ^{-\lambda })\delta ^{-1}\sum_{s=n}^{n+p-1}{\nu
}_{n+p}\nu _{n+p-1}...\nu _{s+2}\nu _{s+1}^{3}\delta ^{-\lambda s}\tilde{%
\varphi}_{s}^{(2)}.
\end{equation}
Then the convergence of the product (\ref{C0}) yields
\begin{equation}
a_{n,p}-1<\left( \prod_{k=n+1}^{\infty }\nu _{k}\right) -1\rightarrow 0,\ \
n\rightarrow +\infty .  \label{liman}
\end{equation}
On the other hand, (\ref{C0}) and (\ref{efi2}) give
\begin{equation}
\left| b_{n,p}\right| <(\theta +1)(1-\delta ^{-\lambda })\delta
^{-1}K_{0}^{7}\left| \tilde{\varphi}_{0}^{(2)}\right| \sum_{s=n}^{\infty
}\delta ^{-(1+\lambda )s}\rightarrow 0,\ \ n\rightarrow +\infty .
\label{limbn}
\end{equation}
The estimates (\ref{e2}) and (\ref{e3}) yield respectively
\begin{equation}
\tilde{\varphi}_{n+p}^{(1)}<a_{n,p}\tilde{\varphi}_{n}^{(1)},\qquad \tilde{%
\varphi}_{n+p}^{(1)}>a_{n,p}\tilde{\varphi}_{n}^{(1)}+b_{n,p}.
\end{equation}
Therefore
\begin{equation}
(a_{n,p}-1)\tilde{\varphi}_{n}^{(1)}+b_{n,p}<\tilde{\varphi}_{n+p}^{(1)}-%
\tilde{\varphi}_{n}^{(1)}<(a_{n,p}-1)\tilde{\varphi}_{n}^{(1)}.  \label{e9}
\end{equation}
Having in mind (\ref{Z8}) and (\ref{e8}), one gets
\begin{equation}
0<\tilde{\varphi}_{n}^{(1)}=\delta ^{\lambda n}{\varphi}_{n}^{\left(
1\right) }<K.  \label{e10}
\end{equation}
Now it suffices to apply the latter estimate together with (\ref{liman}) and (\ref
{limbn}) in (\ref{e9}) and conclude that $\{\tilde{\varphi}_{n}^{(1)}\}$ is a Cauchy sequence.
Thus, for every $\beta <\beta _{\ast }$, there exists $%
\tilde{\varphi}=\tilde{\varphi}(\beta )>0$ such that $\tilde{\varphi}%
_{n}^{(1)}\rightarrow \tilde{\varphi}.$ Now we apply Lemma \ref
{Conver} and obtain the convergence to be proved.
\hfil\vrule height6pt width5pt depth-1pt
\begin{Rk}
When proving the convergence of $\{\tilde{\varphi}_{n}^{(1)}\},$ the limit of
this sequence has been estimated. Namely, we have obtained (see (\ref{e10}))
\begin{equation}
\lim_{n\rightarrow \infty }\tilde{\varphi}_{n}^{(1)}\leq K={\varphi}%
_{0}^{(1)}\prod_{n=1}^{\infty }\nu _{n}=\prod_{n=1}^{\infty }\frac{\hat{%
\varphi}_{0}^{(1)}}{1-(1-\delta ^{-\lambda }){\varphi}_{n-1}^{(1)}}.
\end{equation}
This bound is achieved for $f_{0} (z)
= C \exp (\alpha z)$ (in this case we may calculate $\tilde{g}_n $ explicitly,
see (\ref{G3})). It is quite likely that this bound is achieved also in the
general case, but to prove this conjecture we would need more sophisticated
estimates than (\ref{e3}) or (\ref{e6}).
\end{Rk}

\end{document}